\documentclass[a4paper,reqno,11pt]{amsart}

\usepackage[foot]{amsaddr}
\usepackage{amsmath,amsthm,amssymb,dsfont,graphicx,xspace,epsfig}
\usepackage[dvipsnames]{xcolor}
\usepackage[T1]{fontenc}
\usepackage{thm-restate}
\usepackage[hidelinks]{hyperref}
\usepackage{bm}
\usepackage{color}
\usepackage{comment}
\usepackage{mathrsfs}
\usepackage{cleveref}
\usepackage[shortlabels]{enumitem}
\usepackage{float}
\usepackage{pict2e}
\usepackage{mathtools}
\usepackage{tikz}
\usepackage{subfigure}
\usepackage{xstring,refcount}
\usepackage{yfonts}
\usepackage[longnamesfirst,numbers,sort&compress]{natbib}

\usetikzlibrary{fit, backgrounds, matrix, arrows.meta}
\usetikzlibrary{calc,arrows,positioning}
\usetikzlibrary{decorations}
\usetikzlibrary{decorations.pathmorphing}
\usetikzlibrary{decorations.pathreplacing}
\usetikzlibrary{decorations.shapes}
\usetikzlibrary{decorations.text}
\usetikzlibrary{decorations.markings}
\usetikzlibrary{decorations.fractals}
\usetikzlibrary{decorations.footprints}

\tikzset{vertex/.style = {circle,fill=black,minimum size=6pt, inner sep=0pt, outer sep=3pt}}
\tikzset{mvertex/.style = {circle,fill=black,minimum size=5pt, inner sep=0pt, outer sep=2pt}}
\tikzset{svertex/.style = {circle,fill=black,minimum size=4pt, inner sep=0pt, outer sep=2pt}}
\tikzset{labelledvertex/.style = {circle,fill=none,draw,very thick, inner sep=2pt, outer sep=2pt}}
\tikzset{Xrect/.style = {rectangle,fill=none,draw,very thick, minimum width=1.2cm, minimum height=0.8cm, outer sep=0pt}}

\tikzset{vtarc/.style = {line width=2.2pt,->,> = latex}}
\tikzset{tarc/.style = {thick,->,> = latex}}
\tikzset{sarc/.style = {thin,->,> = latex}}
\tikzset{tdigon/.style = {line width=2.2pt,<->,> = latex}}
\tikzset{digon/.style = {thick,<->,> = latex}}
\tikzset{bigvertex/.style = {shape=circle,draw}}

\tikzset{ledge/.style = {gray}}
\tikzset{tedge/.style = {thick}}
\tikzset{vtedge/.style = {line width=1.3pt}}
\tikzset{stedge/.style = {line width=2.2pt}}

\definecolor{g-blue}{rgb}{0.0, 0.5, 1.0}
\definecolor{g-green}{rgb}{0.3, 0.8, 0.3522}

\newtheorem{theorem}{Theorem}[section]
\newtheorem{lemma}[theorem]{Lemma}
\newtheorem{corollary}[theorem]{Corollary}
\newtheorem{proposition}[theorem]{Proposition}
\newtheorem{conjecture}[theorem]{Conjecture}

\newtheorem{claim}{Claim}
\newcounter{oldthm}{}
\def\@extractthmnum#1.#2{#2}

\theoremstyle{definition}

\newtheorem{problem}[theorem]{Problem}

\newenvironment{proofclaim}[1][Proof of claim]
  {\begin{proof}[#1]}
  {\end{proof}}

\newcommand{\Ra}{\Rightarrow}
\newcommand{\dic}{\vec{\chi}}

\newcommand{\delmax}{\Delta_{\max}}
\newcommand{\delmin}{\Delta_{\min}}
\newcommand{\deltil}{\tilde{\Delta}}
\newcommand{\delplus}{\Delta^{\!+}}

\newcommand{\EE}{\mathbb{E}}
\newcommand{\PP}{\mathbb{P}}

\DeclareMathOperator{\UG}{UG}

\newcommand{\rev}[1]{\reflectbox{\ensuremath{\vec{\reflectbox{\ensuremath{#1}}}}}}
\newcommand{\bid}{\overleftrightarrow}

\newcommand{\bidbistiny}[1]{\overset{\text{\scalebox{.8}{$\bm\leftrightarrow$}}}{#1}}

\newcommand{\bic}{\bidbistiny{\omega}}
\newcommand{\bictiny}{\bidbistiny{\omega}}

\newcommand{\ind}[1]{[#1]}
\newcommand{\edeg}{m}
\newcommand{\edegp}{m^+}
\newcommand{\edegm}{m^-}

\DeclareMathOperator{\IG}{IG}

\renewcommand{\epsilon}{\varepsilon}
\renewcommand{\emptyset}{\varnothing}
\renewcommand{\phi}{\varphi}

\let\leq\leqslant
\let\geq\geqslant

\hypersetup{
    pdftitle={An analogue of Reed's conjecture for digraphs},
    colorlinks,
    linkcolor={red!50!black},
    citecolor={blue!50!black},
    urlcolor={blue!80!black}
}

\title{An analogue of Reed's conjecture for digraphs}

\author[K. Kawarabayashi]{Ken-ichi Kawarabayashi}
\address[K.Kawarabayashi]{National Institute of Informatics, The University of Tokyo, Tokyo, Japan}
\email{k\_keniti@nii.ac.jp}

\author[L. Picasarri-Arrieta]{Lucas Picasarri-Arrieta}
\address[L. Picasarri-Arrieta]{National Institute of Informatics, Tokyo, Japan}
\email{lpicasarr@nii.ac.jp}

\thanks{Research supported by JSPS Kakenhi 26K21777 and JP25K24465 and by JST ASPIRE JPMJAP2302.}
\thanks{A shorter version of this paper has been published in the proceedings of SODA 2025~\cite{kawarabayashiSODA25}.}

\begin{document}

\begin{abstract}
In 1998, Reed conjectured that every graph $G$ satisfies 
$\chi(G) \leq \lceil \frac{\Delta(G)+1+\omega(G)}{2} \rceil$, and proved that,  for some $\epsilon>0$, $\chi(G) \leq \lceil (1-\epsilon)(\Delta(G)+1)+\epsilon\omega(G) \rceil$.

We propose an analogous conjecture for digraphs.
Given a digraph $D$, we denote by $\dic(D)$ the dichromatic number of $D$, which is the minimum number of colours needed to partition $V(D)$ into sets inducing acyclic subdigraphs.
A {\it biclique} is a set of vertices that are
pairwise linked with two arcs in opposite directions, and  $\bic(D)$ denotes the order of a largest biclique of $D$. We also let $\displaystyle\deltil(D) = \max_{v\in V(D)} \sqrt{d^+(v) \cdot d^-(v)}$. 

We conjecture that every digraph $D$ satisfies $\dic(D) \leq \lceil \frac{\deltil(D)+1+\bictiny(D)}{2} \rceil$, which, if true, implies Reed's conjecture.
As a partial result, we prove that there exists $\epsilon >0$ such that every digraph $D$ satisfies $\dic(D) \leq \lceil (1-\epsilon)(\deltil(D)+1)+\epsilon\bic(D) \rceil$. This implies both Reed's result and an independent result of Harutyunyan and Mohar for oriented graphs.  

To obtain this upper bound on $\dic$, we prove that every digraph $D$ whose biclique number is larger than $\frac{2}{3}(\delmax(D)+1)$ admits an acyclic set of vertices intersecting each maximum biclique of $D$, where $\displaystyle\delmax(D) = \max_{v\in V(D)} \max(d^+(v),d^-(v))$. This generalises a result of King.

We finally give a short proof that every oriented graph $D$ with underlying graph $G$ satisfies both $\dic(D) \leq \frac{\sqrt{2}}{2} \deltil(D) + 2$ and $\dic(D) \leq \frac{1}{3}\Delta(G)+2$, improving on results of Golowich and Steiner, respectively.

\end{abstract}

\maketitle

\clearpage

\section{Introduction}

In 1998, Reed~\cite{reedJGT27} posed the following celebrated conjecture, which draws a connection between the maximum degree $\Delta(G)$, the clique number $\omega(G)$, and the chromatic number $\chi(G)$ of a graph~$G$.

\begin{conjecture}[{Reed~\cite[Conj.~1]{reedJGT27}}]
    \label{conj:reed}
    Every graph $G$ satisfies 
    \[
    \chi(G) \leq \left\lceil \frac{\Delta(G) + 1 + \omega(G)}{2} \right\rceil.
    \]
\end{conjecture}

As a partial result towards this conjecture, Reed proved the following.

\begin{theorem}[{Reed~\cite[Cor.~2]{reedJGT27}}]
    \label{thm:reed}
    There exists $\epsilon >0$ such that every graph $G$ satisfies 
    \[ \chi(G) \leq \lceil (1-\epsilon)(\Delta(G) +1) + \epsilon\omega(G)\rceil.\]
\end{theorem}

Observe that Conjecture~\ref{conj:reed} is exactly Theorem~\ref{thm:reed} for $\epsilon = \frac{1}{2}$. 
For graphs with sufficiently large maximum degree, Hurley \textit{et al.}~\cite{hurleyAC22} proved that Theorem~\ref{thm:reed} holds for $\epsilon \leq 0.119$, improving on earlier results obtained successively by Bonamy \textit{et al.}~\cite{bonamyJCTB155} and  by Delcourt and Postle~\cite{delcourtENDM61}.

The purpose of this work is to propose analogues of Conjecture~\ref{conj:reed} and Theorem~\ref{thm:reed} for digraphs. We thus need three digraph parameters that extend $\Delta$, $\omega$, and $\chi$, respectively. By extension, we mean that if $D$ is a symmetric digraph, these digraph parameters on $D$ must coincide with $\Delta$, $\omega$, and $\chi$ on the underlying graph of $D$. One can easily check that this is actually the case for $\deltil$, $\bic$, and $\dic$, defined as follows.

Let $D$ be a digraph. We let $\deltil(D)$ denote the maximum geometric mean of the in- and out-degrees of the vertices of $D$, that is $\deltil(D) = \max \{ \sqrt{d^+(v)\cdot d^-(v)} : v \in V(D) \}$.
A {\it complete digraph} is a digraph whose vertices are
pairwise linked with two arcs in opposite directions.
A \textit{biclique} of $D$ is a set of vertices inducing a complete digraph.
The \textit{biclique number} $\bic(D)$ of $D$ is the order of a largest biclique of $D$.
The notions of dicolouring and dichromatic number were introduced in the late 1970s by Erd\H{o}s and Neumann-Lara~\cite{neumannlaraJCT33,erdosPNCN1979}.
A \textit{$k$-dicolouring} of $D$ is a function $c\colon V(D) \to [k]$ such that $c^{-1}(i)$ induces an acyclic subdigraph of $D$ for each $i \in [k]$. The \textit{dichromatic number} of $D$, denoted by $\dic(D)$, is the smallest $k\in \mathbb{N}$ such that $D$ admits a $k$-dicolouring.
These two concepts were rediscovered by Mohar in the 2000s~\cite{moharJGT43,bokalJGT46} and received a lot of attention since then. It turned out that many classical results on graph colouring can be extended to digraphs with these notions, which shows that dicolouring is somehow the right extension of proper colouring to digraphs. For instance, people generalised the Strong Perfect Graph Theorem~\cite{andresJGT79}, Gy\'arf\'as-Sumner's Conjecture~\cite{aboulkerEJC28,cookEJC30,aboulkerEJC31b}, and
Brooks-type results~\cite{harutyunyanSIDMA25, aboulkerDM346, aboulker2023digraph, goncalvesJGTtoappear}. 

Harutyunyan and Mohar~\cite{harutyunyanEJC18} posed the following analogue of Conjecture~\ref{conj:reed} for oriented graphs, that is, for digraphs with biclique number at most~$1$.

\begin{conjecture}[{Harutyunyan and Mohar~\cite[Conj.~1.5]{harutyunyanEJC18}}]
    \label{conj:harutyunyan}
    Every oriented graph $D$ satisfies
     \[ \dic(D) \leq \left\lceil \frac{\deltil(D)}{2} \right\rceil +1.\]
\end{conjecture}

As a partial result towards their conjecture, they obtained the following.

\begin{theorem}[{Harutyunyan and Mohar~\cite[Thm.~2.1]{harutyunyanEJC18}}]
    \label{thm:harutyunyan}
    There exists $\epsilon>0$ such that every oriented graph $D$ with $\deltil(D)$ large enough satisfies $\dic(D) \leq (1-\epsilon) \deltil(D)$.
\end{theorem}

Theorem~\ref{thm:harutyunyan} was originally proved to hold for $\epsilon=e^{-13}$.
This has been strengthened by Golowich~\cite{golowichDM339}, who obtained that all oriented graphs $D$ satisfy $\dic(D)\leq \sqrt{\frac{2}{3}}\deltil(D) + O(1)$. 
Combining Golowich's idea with a recent result of the second author, in Section~\ref{sec:oriented_graphs} we give a short proof of the following stronger bound, thereby approaching Conjecture~\ref{conj:harutyunyan}.

\begin{restatable}{theorem}{thmorientedgraphs}
    \label{thm:oriented_graphs}
    Let $D$ be an oriented graph with underlying graph $G$. Then  $\dic(D) \leq \frac{1}{3}\Delta(G)+2$, $\dic(D) \leq \frac{2}{3}\delmax(D)+2$, and $\dic(D) \leq \frac{\sqrt{2}}{2}\deltil(D) + 2$.
\end{restatable}

We note that Theorem~\ref{thm:oriented_graphs} also implies a result of Steiner~\cite{steinerDAM287}, stating that if $D$ is a digraph with underlying graph $G$ and digirth $g$, then $\dic(D) \leq \frac{1}{3}\Delta(G) + O\left(\frac{\Delta(G)}{g} + g\right)$.

\medskip

We actually believe that the analogue of Conjecture~\ref{conj:reed} holds not only for oriented graphs but for all digraphs. We thus propose the following conjecture, which, if true, implies both Conjectures~\ref{conj:reed} and~\ref{conj:harutyunyan}.

\begin{conjecture}
    \label{conj:main}
    Every digraph $D$ satisfies 
    \[
    \dic(D) \leq \left\lceil \frac{\deltil(D) + 1 + \bic(D)}{2} \right\rceil.
    \]
\end{conjecture}

The main result of this work is the following, which implies both Theorem~\ref{thm:reed} (when restricted to symmetric digraphs) and Theorem~\ref{thm:harutyunyan} (when restricted to oriented graphs), and provides some support to Conjecture~\ref{conj:main}.

\begin{restatable}{theorem}{mainthm}
    \label{thm:main}
    There exists $\epsilon >0$ such that every digraph $D$ satisfies 
    \[ \dic(D) \leq \lceil (1-\epsilon)(\deltil(D) +1) + \epsilon\bic(D)\rceil. \]
\end{restatable}

We made no attempt to maximize the value of $\epsilon$ in our proof. In particular, applying Theorem~\ref{thm:main} to obtain Theorems~\ref{thm:reed} and~\ref{thm:harutyunyan} yields values of $\epsilon$ that are much smaller than the existing ones.

Our proof of Theorem~\ref{thm:main} is inspired by a short proof of Theorem~\ref{thm:reed} due to King and Reed~\cite{kingJGT81}. Along the way, we generalise a collection of intermediate results to digraphs that we find of independent interest. We refer the reader to the recent book of Stiebitz, Schweser, and Toft~\cite[Chapter~6]{stiebitz2024} for a complete proof of Theorem~\ref{thm:reed} which includes proofs of all intermediate results (on undirected graphs) mentioned later on.

In Section~\ref{sec:sparse_digraphs}, using the probabilistic method, we prove the following upper bound on the dichromatic number of sparse digraphs, which generalises to digraphs a classical result of Molloy and Reed~\cite[Theorem~10.5]{molloyreed}. Given a digraph $D$ with maximum degree $\delmax(D) = \Delta$ and a vertex $v$ of $D$, we denote by $\edegp(v)$ and $\edegm(v)$ the number of arcs of $D\ind{N^+(v)}$ and $D\ind{N^-(v)}$ respectively. We say that $D$ is \textit{$B$-sparse} if, for every vertex $v\in V(D)$, we have $\min(\edegp(v),\edegm(v)) \leq \Delta(\Delta-1) - B$.

\begin{restatable}{theorem}{thmsparsedigraphs}
    \label{thm:sparse_digraphs}
    There exists an integer $\Delta_0 \geq 1$ such that the following holds.
    For every digraph $D$ with maximum degree $\delmax(D) = \Delta \geq \Delta_0$, if $D$ is $B$-sparse with $B > \Delta \log^3 \Delta$, then 
    \[ \dic(D) \leq (\Delta+1) - \frac{B}{4e^7\Delta}. \]
\end{restatable}

In Section~\ref{sec:transversal}, we give a sufficient condition for a digraph to admit an acyclic set of vertices intersecting each of its maximum bicliques, by proving the following generalisation of a result of King~\cite{kingJGT67}.
\begin{theorem}
    Let $D$ be a digraph with $\bic(D) > \frac{2}{3}(\delmax(D)+1)$. Then $D$ has an acyclic set of vertices $I$ such that $\bic(D-I) = \bic(D) -1$.
\end{theorem}
This implies that, in a minimum counterexample to Theorem~\ref{thm:main}, there must be a significant gap between the biclique number and the maximum degree.
We actually prove a more general result, as we characterise exactly the digraphs $D$ with $\bic(D) \geq \frac{2}{3}(\delmax(D)+1)$ that do not admit such an acyclic set of vertices. We thus extend a result due to Christofides {\it et al.}~\cite{christofidesJGT73}. We also discuss the consequences of our result on list-dicolouring of digraphs, and show a generalisation of a result due to Haxell~\cite{haxellCPC10} that we find of independent interest.

In Section~\ref{sec:dense_neighbourhoods}, we deal with vertices having a dense neighbourhood in a counterexample to Theorem~\ref{thm:main}. We extend a result of King and Reed~\cite{kingJGT81} by proving the following.

\begin{restatable}{theorem}{thmdensevertex}
    \label{thm:handle_dense_vertex}
    Let $a,\epsilon>0$ be real numbers with $\epsilon \leq \frac{1}{6} - 4\sqrt{a}$. There exists $\Delta(a)\geq 1$ such that the following holds. Let $D$ be a digraph with maximum degree $\delmax(D) = \Delta \geq \Delta(a)$ satisfying $\bic(D) \leq \frac{2}{3}(\Delta+1)$ and let $v\in V(D)$ be such that 
    \[
    \max(\edegm(v),\edegp(v)) > (1-a)\Delta(\Delta-1).
    \]
    Then $\dic(D) \leq \max\left( \dic(D-v), (1-\epsilon)(\Delta+1) \right)$.
\end{restatable}

We put the pieces together and prove Theorem~\ref{thm:main} in Section~\ref{sec:proof_main}. In particular, we use a trick introduced by Harutyunyan and Mohar~\cite{harutyunyanEJC18} to obtain the result with $\deltil$ instead of $\delmax$.
We finally propose some further research directions in Section~\ref{sec:conclusion}. In particular, we discuss possible analogues of Conjecture~\ref{conj:main} with $\deltil(D)$ being replaced by the smaller parameter $\displaystyle\delmin(D) = \max_{v\in V(D)} \min(d^-(v),d^+(v))$.

\section{Preliminaries}
\label{sec:preliminaries}

\subsection{General notation on digraphs}
Our notation on digraphs follows~\cite{bang2009}.
 A \textit{digon} is a pair of arcs in opposite directions between the same vertices. A \textit{simple arc} is an arc that is not in a digon.
A digraph is \textit{symmetric} if it contains no simple arc.
An \textit{oriented graph} is a digraph with no digon.
Given an undirected graph $G$, the digraph $\bid{G}$ is the symmetric digraph obtained from $G$ by replacing every edge with a digon. For an integer $n$, $\bid{K_n}$ is called the \textit{complete digraph} on $n$ vertices.

Let $D$ be a digraph.
 The \textit{underlying graph} of $D$, denoted by $\UG(D)$, is the undirected graph with vertex-set $V(D)$ in which $uv$ is an edge if and only if $uv$ or $vu$ is an arc of $D$.
 The \textit{symmetric part} of $D$, denoted by $S(D)$, is the undirected graph with vertex-set $V(D)$ in which $uv$ is an edge if and only if $\{uv,vu\}$ is a digon of $D$. A \textit{matching} of $D$ is a set of pairwise disjoint arcs of $D$.
 The \textit{complement} $\bar{D}$ of $D$ is the digraph with vertex-set $V(\bar{D}) = V(D)$ and arc-set $A(\bar{D})=\{uv : u\neq v \text{~and~} uv\notin A(D)\}$.

Let $v$ be a vertex of $D$. The {\it out-degree} (resp. {\it in-degree}) of $v$, denoted by $d^+(v)$ (resp. $d^-(v)$), is the number of arcs leaving (resp. entering) $v$. 
We define the \textit{maximum degree} of $v$ as $d_{\max}(v) = \max\{d^+(v), d^-(v)\}$ and its \textit{minimum degree} as $d_{\min}(v) = \min\{d^+(v), d^-(v)\}$. We can then define the corresponding maximum degrees of $D$: 
\[
    \delmax(D) = \max_{v\in V(D)} (d_{\max}(v)) \hspace{0.8cm} \text{and} \hspace{0.8cm} \delmin(D) = \max_{v\in V(D)} (d_{\min}(v)).
\]
The \textit{maximum geometric mean} of $D$ is defined as $\deltil(D) = \max_{v\in V(D)} \sqrt{d^-(v) \cdot d^+(v)}$. By definition, observe that we have $\delmin(D) \leq \deltil(D) \leq \delmax(D)$. We say that $D$ is {\it $\Delta$-diregular} if, for every vertex $v\in V(D)$, $d^-(v) = d^+(v) = \Delta$. Note that $D$ being diregular implies $\delmin(D) = \delmax(D)$, but the converse is not true.

Let $X,Y$ be two disjoint sets of vertices of $D$. We denote by $A(X,Y)$ the set of arcs $uv\in A(D)$ such that $u\in X$ and $v\in Y$.

A \textit{list assignment} $L$ of $D$ is a map that associates a list of colours to every vertex $v$ of $D$. It is a \textit{$k$-list assignment} if $|L(v)| \geq k$ holds for every vertex $v$ of $D$.
Let $L$ be a list assignment. An \textit{$L$-dicolouring} of $D$ is a dicolouring $\alpha$ of $D$ such that, for every vertex $v$ of $D$, $\alpha(v) \in L(v)$.
We say that $D$ is \textit{$L$-dicolourable} if $D$ admits an $L$-dicolouring.
Given an integer $k$, we say that $D$ is \textit{$k$-dichoosable} if, for every $k$-list assignment $L$ of $D$, $D$ is $L$-dicolourable.
\subsection{Probabilistic tools}

Our proof of Theorem~\ref{thm:main} contains some probabilistic arguments. We briefly recall here two classical results used later. 
The reader unfamiliar with the probabilistic method is referred to~\cite{alon2008} and~\cite{molloyreed}.
We first need the symmetric version of the Lov\'asz Local Lemma, due to Erd\H{o}s and Lov\'asz~\cite{erdosIFS10}.

\begin{lemma}[Lov\'asz Local Lemma]
\label{lemma:lll}
    Let $A_1,A_2,\dots,A_n$ be events in an arbitrary probability space. Suppose that each event $A_i$ is mutually independent of all the other events but at most $d$ of them, and that $\PP(A_i) \leq p$ for all $1\leq i \leq n$. If 
    \[ ep(d+1) \leq 1,\]
    then $\PP\left(\bigwedge_{i=1}^n \overline{A_i}\right)>0$.
\end{lemma}

The second result we need is the celebrated concentration bound due to Talagrand~\cite{talagrand1995}. We use it in the following form (for a proof that it is actually derived from Talagrand's original inequality, see the appendix of~\cite{harutyunyanArxiv25}).

\begin{lemma}[{Talagrand's Inequality}]
    \label{lemma:talagrand}
    Let $X$ be a random variable valued in $\mathbb{N}$, determined by $n$ independent trials and satisfying the following for some integers $c,r \geq 1$:
    \begin{enumerate}
        \item changing the outcome of any one trial can affect $X$ by at most $c$, and
        \item for every $s\in \mathbb{N}$, if $X\geq s$ then there is a set of at most $rs$ trials whose outcomes certify that $X\geq s$.
    \end{enumerate}
    Then
    \[ 
    \PP\left(|X- \EE(X)| > t \right) \leq 4\exp\left(\frac{-t^2}{32c^2 r(\EE(X) + t)} \right)
    \]
    for any real number $t>126c\sqrt{r\EE(X)} + 344 c^2 r$.
\end{lemma}

\section{A short proof for oriented graphs}
\label{sec:oriented_graphs}

In this section we give a short proof of Theorem~\ref{thm:oriented_graphs}, which we first recall here for convenience.

\thmorientedgraphs*

Our proof relies on a combination of the two following known results.

\begin{theorem}[{Lov\'asz~\cite[Thm.~1]{lovaszSSMH1}}]
    \label{thm:lovasz}
    Let $G$ be a graph and $d_1\geq \ldots \geq d_p\geq 1$ be a sequence of integers such that $\sum_{i=1}^pd_i \geq \Delta(G)-p+1$. Then $G$ admits a partition $(V_1,\ldots,V_p)$ of its vertex-set such that, for every $i\in [p]$, $\Delta(G[V_i]) \leq d_i$.
\end{theorem}

\begin{theorem}[{Picasarri-Arrieta~\cite[Cor.~8]{picasarriJGT106}}]
    \label{thm:brooks_delmin}
    Every oriented graph $D$ satisfies 
    \[
        \dic(D) \leq \max(2,\delmin(D)).
    \]
\end{theorem}

\begin{proof}[Proof of Theorem~\ref{thm:oriented_graphs}]
    Let $D$ be any oriented graph with underlying graph $G$. By Theorem~\ref{thm:lovasz} applied with $p=\lceil(\Delta(G)+1)/6 \rceil$ and $d_1=\ldots=d_p=5$, $G$ admits a partition $(V_1,\dots,V_p)$ of its vertex-set such that $\Delta(G[V_i]) \leq 5$.

    For every $i\in [p]$, we obtain that $D[V_i]$ is an orientation of a graph with maximum degree at most $5$. In particular, it follows that $\delmin(D[V_i]) \leq 2$. By Theorem~\ref{thm:brooks_delmin}, we thus have $\dic(D[V_i]) \leq 2$. We thus obtain
    \begin{equation}
        \dic(D) \leq \sum_{i=1}^p \dic(D[V_i]) \leq 2p \leq  \frac{1}{3}\Delta(G)+2.
        \label{eq:Delta/3}
    \end{equation}
    It is straightforward that $\Delta(G) \leq 2\delmax(D)$, and it follows that $\dic(D) \leq \frac{2}{3}\delmax(D)+2$.

    With~\eqref{eq:Delta/3} in hand, we now prove the last inequality of the statement. Assume for a contradiction that there exist a real number $\deltil$ and an oriented graph $D$ with $\deltil(D) \leq \deltil$ and  $\dic(D) > \left\lfloor \frac{\sqrt{2}}{2}\deltil+2 \right\rfloor = k$. Among all such oriented graphs $D$, we choose one 
    with minimum order. This implies that all vertices $v\in V(D)$ satisfy
    \[
    \min(d^-(v),d^+(v)) \geq k\geq \frac{\sqrt{2}}{2}\deltil,
    \]
    for otherwise a $k$-dicolouring of $D-v$ (which exists by minimality of $D$) can be extended to $D$ by choosing for $v$ a colour from $[k]$ that is not appearing in its in- or out-neighbourhood.
    
    Let $v$ be any vertex for which $d^-(v)+d^+(v)$ is maximised, that is $d^-(v)+d^+(v) = \Delta(G)$, where $G$ is the underlying graph of $D$.
    We assume by directional duality that $d^-(v) \geq d^+(v)$, and in particular $d^+(v) \leq \deltil$. By definition of $\deltil$, we have $d^-(v) \leq \frac{\deltil^2}{d^+(v)}$.

    Let $f\colon x \mapsto x+\frac{\deltil^2}{x}$. Observe that, over the interval $[\frac{\sqrt{2}}{2}\deltil, \deltil]$, $f$ is maximised for $x=\frac{\sqrt{2}}{2}\deltil$.
    This shows 
    \[
        \Delta(G) = d^-(v) + d^+(v) \leq d^+(v) + \frac{\deltil^2}{d^+(v)} \leq \frac{3\sqrt{2}}{2} \deltil.
    \]
    It follows from~\eqref{eq:Delta/3} and the inequality above that
    \[
        \dic(D) \leq \frac{\sqrt{2}}{2}\deltil + 2.
    \]
    Since $\dic(D)$ is an integer, it follows that $\dic(D)\leq k$, a contradiction. 
\end{proof}

\section{Colouring sparse digraphs}
\label{sec:sparse_digraphs}

In view of proving Theorem~\ref{thm:main}, our first goal is to generalise the following result due to Molloy and Reed~\cite{molloyreed} on the chromatic number of sparse graphs. For a vertex $v$ of a graph $G$, we denote by $\edeg(v)$ the number of edges of $G\ind{N(v)}$. We say that $G$ is \textit{$B$-sparse} if $\edeg(v) \leq \binom{\Delta(G)}{2} - B$ for every $v\in V(G)$.

\begin{theorem}[{Molloy and Reed~\cite[Thm.~10.5]{molloyreed}}]
    \label{thm:sparse_graphs}
    Let $\Delta$ be an integer large enough and $G$ be a graph with maximum degree $\Delta$. If $G$ is $B$-sparse with $B>\Delta \log^3\Delta$, then 
    \[\chi(G) \leq (\Delta+1) - \frac{B}{e^6\Delta}.\] 
\end{theorem}

The goal of this section is to prove Theorem~\ref{thm:sparse_digraphs}, which is a directed analogue of Theorem~\ref{thm:sparse_graphs}. We first recall it here for convenience.

\thmsparsedigraphs*

We will actually use Theorem~\ref{thm:sparse_digraphs} in the following form in Section~\ref{sec:proof_main}.

\begin{corollary}
    \label{cor:sparse_digraphs}
    There exists $\Delta_1 \geq 2$ such that the following holds. Let $D$ be a digraph with maximum degree $\delmax(D) = \Delta \geq \Delta_1$ and $a\in \mathbb{R}$ be such that $a>\frac{\log^3\Delta}{\Delta-1}$, and for every vertex $v\in V(D)$, $\min(\edegp(v),\edegm(v))\leq (1-a)\Delta(\Delta-1)$. Then
    \[ \dic(D) \leq  \left(1-\frac{a}{5e^7}\right)\left(\Delta+1\right).\]
\end{corollary}
\begin{proof}[Proof of Corollary~\ref{cor:sparse_digraphs}]
    By assumption $D$ is $B_a$-sparse for $B_a= a\Delta(\Delta-1) > \Delta \log^3\Delta$. Hence, by taking $\Delta_1 = \max( \Delta_0,9)$, we obtain by Theorem~\ref{thm:sparse_digraphs} that 
    \[ \dic(D) \leq  \left(1-\frac{a}{4e^7}\right)\left(\Delta+1\right) + \frac{a}{2e^7} \leq \left(1-\frac{a}{5e^7}\right)\left(\Delta+1\right),
    \]
    as desired.
\end{proof}

The remainder of this section is devoted to the proof of Theorem~\ref{thm:sparse_digraphs}.
We first briefly justify that we only have to prove it for $\Delta$-diregular digraphs.

\begin{lemma}
    \label{lemma:diregular_sufficient}
    Let $D$ be a $B$-sparse digraph with maximum degree $\delmax(D) = \Delta$. There exists a $B$-sparse $\Delta$-diregular digraph $H$ such that $\dic(D) \leq \dic(H)$.
\end{lemma}
\begin{proof}
    We repeat the following process until $D$ is $\Delta$-diregular. Add to $D$ a disjoint copy of $\rev{D}$, the digraph obtained from $D$ by reversing every arc. Then for every vertex $v\in V(D)$ such that $d^+_D(v) < \Delta$, we add the arc $vv'$, where $v'$ is the copy of $v$ in $\rev{D}$. Similarly, for every vertex $v$ such that $d^-_D(v) < \Delta$, we add the arc $v'v$. 
    Let $D'$ be the obtained digraph. Note that $\dic(D) \leq \dic(D')$ trivially holds as $D$ is a subdigraph of $D'$.
    Observe also that, for every vertex $v\in V(D)$, we have
    $m^+_D(v) = m^+_{D'}(v) = m^-_{D'}(v')$ and $m^-_D(v) = m^-_{D'}(v) = m^+_{D'}(v')$.
    Hence, as $D$ is $B$-sparse, so is $D'$. 
    
    Observe finally that, for every vertex $v\in V(D)$, $d^+_{D'}(v) = d^-_{D'}(v') = \min(d^+_D(v) + 1,\Delta)$ and $d^-_{D'}(v) = d^+_{D'}(v') = \min(d^-_D(v)+1,\Delta)$.
    Hence, after repeating this process at most $\Delta$ times, we eventually find the desired digraph $H$.
\end{proof}

Let $k$ be a positive integer and $\ell$ be a positive real number.
A \textit{partial $(k,\ell)$-dicolouring} $\alpha$ of a digraph $D$ is a $k$-dicolouring of an induced subdigraph of $D$ such that, for every vertex $v$, at least $\ell$ colours appear twice in $N^-(v)$ or in $N^+(v)$. Formally, we must have
\begin{align*}
    &\left|\{ c\in [k] : \exists u,w \in N^-(v), u\neq w, \alpha(u) = \alpha(w) = c\}\right| \geq \ell\\
    \text{or~~~~} &\left|\{ c\in [k] : \exists u,w \in N^+(v), u\neq w, \alpha(u) = \alpha(w) = c\}\right| \geq \ell.
\end{align*} 

\begin{lemma}
    \label{lemma:partial_dicolouring}
    Let $D$ be a $\Delta$-diregular digraph, $k$ be a positive integer, and $\ell\geq0$ be a real number such that $k \leq \Delta+1-\ell$. If $D$ admits a partial $(k,\ell)$-dicolouring, then $\dic(D) \leq \Delta+1 - \ell$.
\end{lemma}
\begin{proof}
    We assume that $\ell$ is an integer, for otherwise we may prove the result for $\lceil\ell\rceil$ instead.
    Assume that $D$ admits a partial $(k,\ell)$-dicolouring.
    We let $\alpha$ be a partial $(\Delta+1-\ell,\ell)$-dicolouring of $D$ for which the number of uncoloured vertices is minimum, the existence of which is guaranteed as $k\leq \Delta+1-\ell$. We claim that $\alpha$ is a dicolouring of $D$, hence implying $\dic(D) \leq \Delta+1 - \ell$. Assume for a contradiction that there is a vertex $v$ which is not coloured by $\alpha$.
    By definition of a partial $(\Delta+1-\ell,\ell)$-dicolouring, the number of colours appearing in both $N^+(v)$ and $N^-(v)$ is at most $\Delta-\ell$.
    Hence we can choose for $v$ a colour of $[\Delta+1-\ell]$ that is not appearing in $N^-(v)$ or in $N^+(v)$. Colouring $v$ with this colour extends $\alpha$ into a partial $(\Delta+1-\ell,\ell)$-dicolouring of $D$ with fewer uncoloured vertices, a contradiction to the choice of $\alpha$.
\end{proof}

We are now ready to prove Theorem~\ref{thm:sparse_digraphs}. 
\begin{proof}[Proof of Theorem~\ref{thm:sparse_digraphs}]
    We do not give the explicit value of $\Delta_0$, we simply assume in what follows that it is arbitrarily large. We assume that $D$ is $\Delta$-diregular, for otherwise we may prove the result for the $B$-sparse $\Delta$-diregular digraph $H$ obtained by Lemma~\ref{lemma:diregular_sufficient}, and conclude that $\dic(D) \leq \dic(H) \leq (\Delta+1) - \frac{B}{4e^7\Delta}$.
    
    We use the probabilistic method. 
    Let us fix $k=\left \lfloor \frac{\Delta}{2} \right \rfloor$. We randomly colour every vertex of $D$ with a colour of $[k]$. Then we simultaneously uncolour every vertex $v$ that has at least one in-neighbour and at least one out-neighbour sharing its colour.
    Our goal is to show that with positive probability, the obtained colouring $\alpha$ is a partial $(k,\ell)$-dicolouring of $D$, where $\ell = \frac{B}{4e^7 \Delta}$. Observe that we have $k\leq \Delta+1-\ell$ as $B\leq \Delta(\Delta-1)$ and $\Delta$ is large enough. Hence, if $\alpha$ is indeed a partial $(k,\ell)$-dicolouring of $D$, the result follows by Lemma~\ref{lemma:partial_dicolouring}. Note that $\alpha$ is necessarily a dicolouring of an induced subdigraph of $D$, for otherwise $D$, partially coloured with $\alpha$, contains a monochromatic directed cycle, the vertices of which should have been uncoloured in the process. It remains to show that, for every vertex of $D$, at least $\ell$ colours are repeated in its in-neighbourhood or in its out-neighbourhood. 
    
    For each vertex $v$, we let $N_v$ be $N^+(v)$ if $\edegp(v) \leq \edegm(v)$ and $N^-(v)$ otherwise.
    We also let $B_v$ be the number of missing arcs in $D\ind{N_v}$, that is $B_v = \Delta(\Delta-1) - \min(\edegm(v), \edegp(v))$. Note that $B_v \geq B > \Delta \log^3 \Delta$.
    Let $X_v$ be the random variable corresponding to the number of colours that are assigned to at least two vertices in $N_v$ and are retained by all these vertices. We also define $Y_v$ as the number of colours assigned to at least two vertices in $N_v$ that are not linked by a digon, but not necessarily retained. Finally, we define $Z_v$ as the number of colours assigned to at least two vertices in $N_v$ that are not linked by a digon and one of these vertices does not retain its colour. We thus have $X_v = Y_v - Z_v$.

    We let $A_v$ be the event that $X_v < \ell$. Our goal is to show that, with positive probability, none of the events $A_v$ occur. To do so, we prove that, for each vertex $v$, the probability that $A_v$ occurs is small enough and then conclude by the Lov{\'a}sz Local Lemma (Lemma~\ref{lemma:lll}). 
    
    \begin{claim}
        \label{claim:large_exp}
        For every vertex $v$, we have $\EE(X_v) \geq \frac{B_v}{2e^7\Delta}$.
    \end{claim}
    \begin{proofclaim}
        For every pair of vertices $\{u,w\} \subseteq N_v$ and every colour $c\in [k]$, we let $A_{\{u,w\}}^c$ be the event that $\{u,w\}$ are exactly the vertices coloured $c$ in $N_v$ and they both retain their colour. We let $X_{\{u,w\}}^c$ be the binary random variable equal to $1$ if $A_{\{u,w\}}^c$ occurs and $0$ otherwise. We clearly have $X_v \geq \sum_{\{u,w\},c} X_{\{u,w\}}^c$. 
        
        Assume $u$ and $w$ are not linked by a digon. If $uw$ is a simple arc of $D$, then $A_{\{u,w\}}^c$ occurs if and only if both $u$ and $w$ are coloured $c$ and no vertex in $N = (N_v \cup N^-(u) \cup N^+(w)) \setminus \{u,w\}$ receives colour $c$. If the arc $uw$ is missing, then $A_{\{u,w\}}^c$ occurs in particular if both $u$ and $w$ are coloured $c$ and no vertex in $N = (N_v \cup N^+(u) \cup N^-(w))\setminus \{u,w\}$ receives colour $c$. In both cases, as $k=\left\lfloor \frac{\Delta}{2} \right\rfloor$, we have $|N| \leq 3\Delta \leq 7k$, which implies
        \[ 
        \PP(A_{\{u,w\}}^c) \geq \left(\frac{1}{k}\right)^2 \left(1 - \frac{1}{k}\right)^{7k}. 
        \]
        By definition of $B_v$ we have at least $\frac{B_v}{2}$ distinct pairs of vertices in $N_v$ that are not linked by a digon. We also have $k$ choices for the colour $c$. By linearity of expectation, we thus obtain
        \[ \EE(X_v) \geq  k \frac{B_v}{2} \left(\frac{1}{k}\right)^2 \left(1-\frac{1}{k}\right)^{7k} \geq \frac{B_v}{\Delta} \left(1-\frac{1}{k}\right)^{7k} \geq \frac{B_v}{2e^7\Delta},\]
        where in the last inequality we used that $k=\left\lfloor \frac{\Delta}{2} \right\rfloor$ is large enough and that 
        \[
        {\displaystyle \lim_{k \to +\infty}} \left(1-\frac{1}{k}\right)^{7k} = \frac{1}{e^7} > \frac{1}{2e^7}. \qedhere
        \]
    \end{proofclaim}
    
    \begin{claim}
        \label{claim:small_exp}
        For every vertex $v$, we have $\EE(X_v) \leq \EE(Y_v) \leq \frac{3B_v}{\Delta}$.
    \end{claim}
    \begin{proofclaim}
        Since $X_v = Y_v - Z_v$, it is clear that $\EE(X_v) \leq \EE(Y_v)$.
        Given a pair of vertices $\{u,w\} \subseteq N_v$ that are not linked by a digon and a colour $c\in [k]$, we let $Y_{\{u,w\}}^c$ be the binary random variable equal to $1$ if both $u$ and $w$ are coloured $c$ and $0$ otherwise.
        We have $Y_v \leq \sum_{\{u,w\},c} Y_{\{u,w\}}^c$ and $\PP(Y_{\{u,w\}}^c = 1) = \frac{1}{k^2}$. Since we have at most $B_v$ choices for the pair $\{u,w\}$ and exactly $k$ for the colour $c$, by linearity of the expectation we obtain
        \[
            \EE(X_v) \leq \EE(Y_v) \leq kB_v\left(\frac{1}{k}\right)^2 \leq \frac{3B_v}{\Delta},
        \]
        where in the last inequality we used $k = \left\lfloor \frac{\Delta}{2} \right\rfloor \geq \frac{\Delta}{3}$.
    \end{proofclaim}

    \begin{claim}
        \label{claim:small_pr}
        For every vertex $v$, $\PP\left(|X_v - \EE(X_v)| > \log \Delta \sqrt{\EE(X_v)}\right) < \frac{1}{\Delta^5}$.
    \end{claim}
    \begin{proofclaim}
        Let us set $t = \frac{1}{2} \log \Delta \sqrt{\EE(X_v)}$.
        Recall that $X_v = Y_v - Z_v$, so by linearity of the expectation we have $\EE(X_v) = \EE(Y_v) - \EE(Z_v)$, and $|X_v - \EE(X_v)| \leq |Y_v - \EE(Y_v)| + |Z_v - \EE(Z_v)|$. In particular, this implies 
        \begin{equation}
            \PP(|X_v - \EE(X_v)| > 2t) \leq \PP(|Y_v - \EE(Y_v)| > t) + \PP(|Z_v - \EE(Z_v)| > t).\label{eq:XvleqYvplusZv}
        \end{equation}
        We use Talagrand's inequality (Lemma~\ref{lemma:talagrand}) to bound both $\PP(|Y_v - \EE(Y_v)| > t)$ and $\PP(|Z_v - \EE(Z_v)| > t)$. Observe that changing the colour of exactly one vertex affects $Y_v$ by at most $1$. Also if $Y_v\geq s$ then there exists a set of $2s$ vertices whose colours certify $Y_v \geq s$: for each of the $s$ colours we may choose two vertices not linked by a digon sharing this colour. 
        Combining Claims~\ref{claim:large_exp} and~\ref{claim:small_exp} we obtain $\EE(Y_v) \leq 6e^7 \EE(X_v)$. Hence, $\Delta$ being sufficiently large, we have $t > 126 \sqrt{2\EE(Y_v)} + 688$, which allows us to apply Talagrand's inequality (Lemma~\ref{lemma:talagrand}) to $Y_v$ and obtain
        \begin{equation}
            \PP(|Y_v - \EE(Y_v)| > t) \leq 4\exp\left(-\frac{t^2}{64(\EE(Y_v) + t)}\right).\label{eq:PrYv}
        \end{equation}
        We now consider $Z_v$. Changing the colour of exactly one vertex affects $Z_v$ by at most $1$. Also if $Z_v \geq s$ then there exists a set of $4s$ vertices whose colours certify $Z_v \geq s$: for each of the $s$ colours we may choose four vertices $a,b,b^-,b^+$ sharing this colour, $a\neq b$, such that $a$ and $b$ are not linked by a digon in $D$, and $b^-b, bb^+ \in A(D)$ ($a, b^-, b^+$ are not necessarily distinct). 
        Observe that $Z_v \leq Y_v$, so in particular $\EE(Z_v) \leq \EE(Y_v) \leq 6e^7 \EE(X_v)$ by Claims~\ref{claim:large_exp} and~\ref{claim:small_exp}. Hence, since $\Delta$ is sufficiently large, we have $t > 126\sqrt{4\EE(Z_v)} + 1376$, which allows us to apply Talagrand's inequality (Lemma~\ref{lemma:talagrand}) to $Z_v$ and obtain 
        \begin{equation}
            \PP(|Z_v - \EE(Z_v)| > t) \leq 4\exp\left(-\frac{t^2}{128(\EE(Z_v) + t)}\right).\label{eq:PrZv}
        \end{equation}
        Combining~\eqref{eq:XvleqYvplusZv},~\eqref{eq:PrYv}, and~\eqref{eq:PrZv} together with $\EE(Z_v) \leq \EE(Y_v) \leq 6e^7\EE(X_v)$ we obtain
        \begin{align*}
            \PP(|X_v - \EE(X_v)| > 2t) &\leq 8\exp\left(-\frac{t^2}{768e^7\EE(X_v) + 128t}\right)
            = 8\exp\left(-\frac{\log^2\Delta }{3072e^7 + 256 \frac{\log \Delta}{\sqrt{\EE(X_v)}}}\right)\\
            &\leq 8\exp\left(-\frac{\log^2\Delta }{3072e^7 + 256\sqrt{2e^7} \frac{1}{\sqrt{\log \Delta}}}\right),
        \end{align*}
        where in the last inequality we used $\EE(X_v) \geq \frac{B_v}{2e^7\Delta}$ and $B_v \geq \Delta \log^3 \Delta$. As $\Delta$ is sufficiently large, we obtain
        \[
            \PP\left(|X_v - \EE(X_v)| > \log \Delta \sqrt{\EE(X_v)}\right) < \frac{1}{\Delta^5},
        \]
        which concludes the proof of the claim.
    \end{proofclaim}

    By Claim~\ref{claim:large_exp} together with the fact that $B_v \geq B \geq \Delta \log^3 \Delta$, we have 
    \[
        \EE(X_v) - \log \Delta \sqrt{\EE(X_v)} \geq \frac{B_v}{2e^7\Delta}\left(1-\log \Delta \sqrt{\frac{2e^7\Delta}{B_v}}\right) \geq \frac{B}{2e^7\Delta}\left(1-\sqrt{\frac{2e^7}{\log \Delta}}\right).
    \]
    By choice of $\ell$, we thus have $\EE(X_v) - \log \Delta \sqrt{\EE(X_v)} \geq 2\ell \left(1-\sqrt{\frac{2e^7}{\log \Delta}}\right)$, which is larger than $\ell$ for $\Delta$ large enough.
    Hence, we obtain 
    \begin{align*}
    \PP(A_v) = \PP\left(X_v<\ell\right) &\leq \PP\left(X_v < \EE(X_v) - \log \Delta \sqrt{\EE(X_v)}\right)\\
    &\leq \PP\left(|X_v -\EE(X_v)| >  \log \Delta \sqrt{\EE(X_v)}\right).
    \end{align*}
    Therefore, by Claim~\ref{claim:small_pr}, we have $\PP(A_v) < \frac{1}{\Delta^5}$. Observe finally that the value of $X_v$ depends only on the colour of vertices at distance at most $2$ of $v$ in $\UG(D)$. Hence, each event $A_v$ is mutually independent of all other events $A_u$ except at most $(2\Delta)^4=d$ of them. For every vertex $v$, we have $e\PP(A_v)(d+1) \leq e\frac{16\Delta^4+1}{\Delta^5}$, which is smaller than $1$ if $\Delta$ is large enough. Hence, by the Lov{\'a}sz Local Lemma (Lemma~\ref{lemma:lll}), with positive probability, none of the events $A_v$ occur. This concludes the proof. 
\end{proof}

\section{Hitting every maximum biclique with an acyclic set of vertices}
\label{sec:transversal}

Let $G$ be a graph with maximum degree $\Delta(G) \leq 3k-1$. Haxell~\cite{haxellCPC10} proved that, if the vertices of $G$ can be partitioned into $r$ cliques $C_1,\dots,C_r$ of size $\geq 2k$, then $G$ contains an independent set $S$ of size $r$. In particular, $S$ intersects each clique $C_i$. 

Using this result, Rabern~\cite{rabernJGT66} proved that every graph $G$ satisfying $\omega(G) \geq \frac{3}{4}(\Delta(G)+1)$ contains an independent set intersecting every maximum clique of $G$. This was later improved to $\omega(G) > \frac{2}{3}(\Delta(G)+1)$ by King~\cite{kingJGT67}, which is best possible. Christofides, Edwards, and King~\cite{christofidesJGT73} later characterised exactly the graphs $G$ satisfying $\omega(G) = \frac{2}{3}(\Delta(G) + 1)$ that do not admit an independent set intersecting every maximum clique. Combined with King's earlier result, they thus obtained the following. Given two graphs $G$ and $G'$, we denote by $G \circ G'$ the lexicographic product of $G$ by $G'$, which is the graph with vertex-set $V(G) \times V(G')$ and edge-set $\{(u,u')(v,v') : uv\in E(G) \text{~or~} (u=v \text{~and~} u'v'\in E(G'))\}$.
\begin{theorem}[{Christofides {\it et al.}~\cite[Thm.~2]{christofidesJGT73}}]
    \label{thm:christofides}
    Let $G$ be a connected graph with $\omega(G) \geq \frac{2}{3}(\Delta(G)+1)$. Then $G$ has an independent set $I$ such that $\omega(G-I) = \omega(G) -1$, unless $G = C_n\circ K_p$ with $n\geq5$ odd and $p=\frac{1}{2}\omega(G)$.
\end{theorem}

The goal of this section is to prove the following generalisation of  Theorem~\ref{thm:christofides}. 
In the proof of Theorem~\ref{thm:main}, we only use Theorem~\ref{thm:decrease_w} with $\bic(D)$ being larger than $\frac{2}{3}(\Delta+1)$. We find the characterisation of the equality case of independent interest and we thus also include it.

\begin{restatable}{theorem}{thmdecreasew}
    \label{thm:decrease_w}
    Let $D$ be a connected digraph and $\Delta$ be an integer such that $\delmax(D) \leq \Delta$ and $\bic(D) \geq \frac{2}{3}(\Delta+1)$. Then $D$ has an acyclic set of vertices $I$ such that $\bic(D-I) = \bic(D)-1$, unless $D=\bid{C_n\circ K_p}$ with $n\geq 5$ odd and $p=\frac{1}{2}\bic(D)=\frac{1}{3}(\Delta+1)$.
\end{restatable}

We first prove the following result, which generalises to digraphs a more precise version of Haxell's result obtained by King~\cite[Thm.~5]{kingJGT67}.

\begin{theorem}
    \label{thm:acyclic_hittingset}
    Let $k$ be a positive integer and $D$ be a digraph with vertices partitioned into bicliques $C_1,\dots, C_r$. If for every $i\in[r]$ and every $v\in C_i$, $v$ has at most $k$ out-neighbours outside $C_i$ and at most $|C_i|-k$ in-neighbours outside $C_i$, then $D$ contains an acyclic set of size $r$.
\end{theorem}

We note that another independent generalisation of Haxell's result to the directed setting, which involves only the out-degree, has been obtained by Aharoni, Berger, and Kfir~\cite{aharoniJGT59}.

Given a digraph $D$ whose vertices are partitioned into independent sets $V_1,\dots,V_r$, an \textit{acyclic system of representatives} (with respect to $V_1,\dots,V_r$), or \textit{ASR} for conciseness, is an acyclic set of vertices $X$ such that, for every $i\in [r]$, $|V_i\cap X|=1$. If we relax this condition to $|V_i\cap X|\leq 1$, then $X$ is a \textit{partial ASR}.
Note that, in the statement of Theorem~\ref{thm:acyclic_hittingset}, any acyclic set $X$ of size $r$ is necessarily a transversal of $(C_1,\dots,C_r)$, as two vertices of $X$ cannot belong to the same biclique. Therefore, Theorem~\ref{thm:acyclic_hittingset} can be restated in terms of ASR as follows.

\begin{theorem}
    \label{thm:ASR}
    Let $k$ be a positive integer and $D$ be a digraph with vertices partitioned into independent sets $V_1,\dots, V_r$. If every vertex $v\in V_i$ satisfies $d^+(v) \leq k$ and $d^-(v) \leq |V_i|-k$, then $D$ admits an ASR.
\end{theorem}

Following the analogue of Haxell's idea~\cite{haxellCPC10}, Theorem~\ref{thm:ASR} has the following consequence about list-dicolouring that we find of independent interest.

\begin{corollary}
    Let $D$ be a digraph, $L$ be a list assignment of $D$, and $k$ be an integer. If for every vertex $v\in V(D)$ and for every colour $c\in L(v)$, we have 
    \[ |N^+(v)\cap \{u : c\in L(u)\}|\leq k \hspace{0.8cm} \text{and} \hspace{0.8cm} |N^-(v)\cap \{u : c\in L(u)\}|\leq  |L(v)| - k,\] 
    then $D$ is $L$-dicolourable.
\end{corollary}
\begin{proof}
    Let $H$ be the auxiliary digraph where $V(H) = \{ (v,c) : v\in V(D), c\in L(v)\}$ which contains an arc from $(u,c_1)$ to $(v,c_2)$ whenever $c_1=c_2$ and $v\in N^+(u)$. 
    Let us label with $v_1,\dots,v_r$ the vertices of $D$, and for every $i\in [r]$ let $V_i = \{(v_i,c) : c\in L(v_i) \}$. By definition, $V_1,\dots,V_r$ is a partition of $V(H)$ into independent sets. 
    The conditions on $L$ imply that every vertex $v\in V_i$ has out-degree at most $k$ and in-degree at most $|V_i|- k$.
    Hence, by Theorem~\ref{thm:ASR}, $H$ contains an ASR $X=\{(v_1,c_1),\dots,(v_r,c_r)\}$. Let $\alpha$ be the colouring of $D$ which associates to each vertex $v_i$ colour $c_i$. We claim that $\alpha$ is an $L$-dicolouring of $D$. We have $c_i \in L(v_i)$ by definition. Assume for a contradiction that $D$, coloured with $\alpha$, contains a monochromatic directed cycle $u_1,\dots,u_\ell$. Let $c$ be the colour of this cycle, then $(u_1,c),\dots,(u_\ell,c)$ is  directed cycle in $H\ind{X}$, a contradiction.
\end{proof}

\begin{proof}[Proof of Theorem~\ref{thm:ASR}]
    For any subset of indices $I \subseteq [r]$, we denote by $V_I$ the set of vertices $\bigcup_{i\in I}V_i$. Given a set of vertices $X\subseteq V(D)$, we denote by $\iota(X)$ the set of indices $\{i\in [r] : V_i \cap X \neq \emptyset\}$. For a single vertex $v$, we denote by $\iota(v)$ the integer $i$ such that $v\in V_i$. 

    We proceed by induction on $r\geq 1$. The result is trivial when $r=1$; we assume $r\geq 2$. Let us fix any vertex $x_1 \in V_r$, we show as a stronger result that there exists an ASR containing $x_1$. Assume for a contradiction that no such ASR exists.

    Let $I\subseteq [r-1]$ and $X,Y \subseteq (V_I\cup \{x_1\})$. We say that $(I,X,Y)$ is a \textit{good triplet} if each of the following properties holds:
    \begin{enumerate}[label=(\alph*)]
        \item $X$ and $Y$ are two disjoint acyclic sets with $x_1\in X$;\label{cond:ayclic}
        \item $Y$ is an ASR of $D\ind{V_I}$ (with respect to $(V_i)_{i\in I}$);\label{cond:Y_ASR}
        \item every vertex in $Y$ has exactly one in-neighbour in $X$ and every vertex in $X$ has at least one out-neighbour in $Y$, in particular $|X|\leq |Y|$ ; and\label{cond:X_smaller_Y}
        \item every vertex $v\in V_I \cup \{x_1\}$ has an in-neighbour in $X$ or an out-neighbour in $Y$.\label{cond:dominating}
    \end{enumerate}

    We first justify that there exists no good triplet in $D$.
    \begin{claim}
        \label{claim:no_good_triplet}
        $D$ admits no good triplet $(I,X,Y)$.
    \end{claim}
    \begin{proofclaim}
        Assume for a contradiction that $(I,X,Y)$ is a good triplet of $D$. The assumption on the degrees in $D$ implies in particular $d^+(x) \leq k$ for every vertex $x\in X$ and $d^-(y) \leq |V_{\iota(y)}| - k$ for every vertex $y\in Y$. Observe that $|X|\leq |Y| = |I|$ by~\eqref{cond:X_smaller_Y} and~\ref{cond:Y_ASR}. We thus have
        \[  
            \sum_{x\in X}d^+(x) + \sum_{y\in Y}d^-(y) \leq k|I| + \sum_{i\in I}(|V_i|-k) = |V_I|.
        \]
        On the other hand, we know by~\ref{cond:dominating} that $(V_I \cup \{x_1\}) \subseteq (N^+(X) \cup N^-(Y))$. This implies
        \[
            \sum_{x\in X}d^+(x) + \sum_{y\in Y}d^-(y) \geq |N^+(X) \cup N^-(Y)| \geq |V_I|+1,
        \]
        yielding the contradiction.
    \end{proofclaim}

    Our goal is now to show that $D$ admits an infinite sequence $(Y_{i})_{i\geq 1}$ of partial ASR, such that each $Y_i$ is a proper subset of $Y_{i+1}$. This yields a contradiction since such a sequence has length at most $r$.

    Let $R_1$ be an ASR of $D\ind{V_{[r-1]}}$ (with respect to $V_1,\dots,V_{r-1}$) for which the size of $Y_1=Y_1'=N^+(x_1) \cap R_1$ is minimum. The existence of $R_1$ is guaranteed by induction on $r$. Observe that $Y_1 \neq \emptyset$ for otherwise $R_1 \cup \{x_1\}$ is an ASR of $D$ containing $x_1$. We let $X_1$ be the singleton $\{x_1\}$.
    We recursively build $(R_{i+1},X_{i+1},Y_{i+1},Y_{i+1}')$ from $(R_{i},X_{i},Y_{i},Y_{i}')$ in such a way that, for every $i\geq 1$, each of the following properties holds.
    \begin{enumerate}[label=(\arabic*)]
        \item $X_i = \{x_1,\dots,x_i\}$ is acyclic. If $i>1$ then $x_i \in V_{\iota(Y_{i-1})}$.\label{cond:Xi}
        \item $R_i$ is an ASR of $D\ind{V_{[r-1]}}$ such that for every $1\leq j\leq i-1$, $R_i\cap N^+(x_j) = Y_j'$. Subject to this property, $Y_i' = N^+(x_i) \cap R_i$ has minimum size.\label{cond:Ri}
        \item $Y_i = \bigsqcup_{1\leq j \leq i} Y_j'$ and $Y_i' \neq \emptyset$.\label{cond:Yi}
    \end{enumerate}

    Observe that $(R_1,X_1,Y_1,Y_1')$ satisfies~\ref{cond:Xi},~\ref{cond:Ri}, and~\ref{cond:Yi}. Also, the existence of such an infinite sequence $(R_i,X_i,Y_i,Y_i')_{i\geq 1}$ implies the desired infinite sequence $(Y_i)_{i\geq 1}$ of partial ASR, yielding the contradiction. 
    
    Hence we assume that for some $i\geq 1$, there exists $(R_j,X_j,Y_j,Y_j')_{1\leq j \leq i}$ such that, for each $j\in [i]$,~\ref{cond:Xi},~\ref{cond:Ri}, and~\ref{cond:Yi} hold. We construct $(R_{i+1},X_{i+1},Y_{i+1},Y_{i+1}')$ as follows.

    \medskip
    
    Let $x_{i+1}$ be a vertex in $V_{\iota(Y_i)}\setminus (X_i \cup Y_i)$ such that
        \begin{equation}\label{eq:xip1}
            N^-(x_{i+1}) \cap X_i = \emptyset \hspace{0.8cm} \text{and} \hspace{0.8cm} N^+(x_{i+1}) \cap Y_i = \emptyset.\tag{$\star$} 
        \end{equation}
        If such a vertex $x_{i+1}$ exists then $X_{i+1}=\{x_1,\dots,x_{i+1}\}$ is acyclic as $x_{i+1}$ is a source in $D\ind{X_{i+1}}$ and $X_i=\{x_1,\dots,x_i\}$ is acyclic by recurrence, so~\ref{cond:Xi} is maintained.
        Let us justify the existence of $x_{i+1}$. If no such vertex exists, then $X_i$ and $Y_i$ are two disjoint acyclic sets such that $x_1\in X_i$, with $Y_i$ being an ASR of $V_{\iota(Y_i)}$. Combining conditions~\ref{cond:Ri} and~\ref{cond:Yi}, we obtain that $(N^+(x_1),\dots,N^+(x_i))$ partitions $Y_i$ into non-empty parts, so in particular every vertex in $Y_i$ has exactly one in-neighbour in $X_i$ and every vertex in $X_i$ has at least one out-neighbour in $Y_i$. Hence, if no vertex in $V_{\iota(Y_i)}\setminus (X_i \cup Y_i)$ satisfies~\eqref{eq:xip1}, we deduce that $(\iota(Y_i), X_i,Y_i)$ is a good triplet, a contradiction to Claim~\ref{claim:no_good_triplet}.

    \medskip
    
    We let $R_{i+1}$ be an ASR of $D\ind{V_{[r-1]}}$ such that for every $1\leq j\leq i$, $R_{i+1}\cap N^+(x_j) = Y_j'$, and subject to this property, $Y_{i+1}' = N^+(x_{i+1}) \cap R_{i+1}$ has minimum size. The existence of $R_{i+1}$ is guaranteed as $R_i$ is such an ASR. Hence~\ref{cond:Ri} is maintained.

    \medskip
    
    It remains to show that~\ref{cond:Yi} is maintained, that is $Y_{i+1}' \neq \emptyset$ and $Y_{i+1}' \cap Y_i = \emptyset$. First we have $Y_{i+1}' \cap Y_i = \emptyset$ by~\eqref{eq:xip1}. Now assume for a contradiction that $Y_{i+1}' = \emptyset$, which together with~\eqref{eq:xip1} implies that $x_{i+1}$ has no out-neighbour in $R_{i+1}$.

    By recurrence we know that $(Y_1',\dots,Y_i')$ partitions $Y_i$. As $x_{i+1}$ belongs to $V_{\iota(Y_i)}\setminus (X_i \cup Y_i)$ there exists a unique $j\leq i$ such that $\iota(x_{i+1}) \in \iota(Y_j')$. As $Y_j'$ is a partial ASR, there exists a unique $y_j \in Y_j' \cap V_{\iota(x_{i+1})}$. Consider the set $\tilde{R}_j$ defined as
    \[
        \Tilde{R}_j = R_{i+1} \setminus \{y_j\} \cup \{x_{i+1}\}.
    \]
    Observe that $\Tilde{R}_j$ is necessarily acyclic as $x_{i+1}$ has no out-neighbour in $R_{i+1}$ and $R_{i+1}$ is an ASR of $D\ind{V_{[r-1]}}$. Also, since $y_j \in V_{\iota(x_{i+1})}$, we obtain that $\tilde{R}_j$ is an ASR of $D\ind{V_{[r-1]}}$.

    Since $x_{i+1}$ has no in-neighbour in $X_i$ by~\eqref{eq:xip1},
    for every $1\leq \ell < j$, we have $\tilde{R}_j\cap N^+(x_{\ell}) = R_{i+1} \cap N^+(x_{\ell}) = Y_{\ell}'$. We finally have $\tilde{R}_j\cap N^+(x_{j}) = (R_{i+1} \cap N^+(x_j))\setminus \{y_j\} = Y_{j}'\setminus \{y_j\}$, a contradiction to the choice of $R_j$.
    Hence~\ref{cond:Yi} is maintained and the proof follows.
\end{proof}

With Theorem~\ref{thm:acyclic_hittingset} in hand, we now prove Theorem~\ref{thm:decrease_w}, for which we need two additional known results. 
Let $G$ be a (di)graph. We denote by $\mathcal{C}(G)$ the set of all maximum (bi)cliques of $G$. Given a collection $\mathcal{C} \subseteq \mathcal{C}(G)$ of maximum (bi)cliques of $G$, we denote by $\IG(\mathcal{C})$ the intersection graph of $\mathcal{C}$, \textit{i.e.}, the undirected graph with vertex-set $\mathcal{C}$ which contains an edge between $C_1$ and $C_2$ whenever $C_1\cap C_2 \neq \emptyset$. We first need the following result due to Hajnal~\cite{hajnalCJM17} (see also~\cite[Lem.~6.11]{stiebitz2024}).
\begin{lemma}[Hajnal~\cite{hajnalCJM17}]
    \label{lemma:hajnal}
    Let $G$ be a graph and $\mathcal{C}$ be a collection of maximum cliques of $G$, then
    \[
         \left| \bigcap \mathcal{C} \right| + \left| \bigcup \mathcal{C} \right| \geq 2\omega(G).
    \]
\end{lemma}

 The second intermediate result we need is an extension of a result of Kostochka~\cite{kostochkaMDA35} due to Christofides, Edwards, and King~\cite{christofidesJGT73}. The statement we use slightly differs from the original one (see~\cite[Lem.~6.12]{stiebitz2024} for a proof of the result in the following form).

\begin{lemma}[Christofides {\it et al.}~\cite{christofidesJGT73}] 
    \label{lemma:christofides}
    Let $G$ be a graph with $\omega(G) \geq \frac{2}{3}(\Delta(G)+1)$, let $\mathcal{C}\subseteq \mathcal{C}(G)$ such that $\IG(\mathcal{C})$ is connected, let $X=\bigcap \mathcal{C}$ and let $Y=\bigcup \mathcal{C}$. 
    If $X=\emptyset$ then $\omega(G) = \frac{2}{3}(\Delta(G)+1)$ and $Y$ has a partition $\{Q_i : i\in [n]\}$ into $n\geq 4$ cliques of $G$ each of cardinality $\frac{1}{2}\omega(G)$ such that $\mathcal{C}\setminus \{Q_1\cup Q_n\} = \{Q_i\cup Q_{i+1} : i \in [n-1]\}$.
\end{lemma}

We are now ready to prove Theorem~\ref{thm:decrease_w}, which we first restate here for convenience.

\thmdecreasew*
\begin{proof}
    Let $D$ be a minimum counterexample, which means that for some $\Delta\in \mathbb{N}$, $\delmax(D) \leq \Delta$, $\bic(D) \geq \frac{2}{3}(\Delta+1)$, $D$ is distinct from the obstruction $\bid{C_n\circ K_p}$ and $D$ does not contain any acyclic set of vertices $I$ intersecting every maximum biclique of $D$. 
    
    We first show that every vertex $v$ belongs to a maximum biclique of $D$. Assume that this is not the case for some vertex $v$ and let $D'$ be $D-v$. Observe then that $\bic(D') = \bic(D) \geq \frac{2}{3}(\Delta+1)$ and $\Delta \geq \delmax(D')$. By minimality of $D$, the statement holds for each connected component $K$ of $D'$ with $\bic(K) = \bic(D')$. 
    Hence, either $D'$ has a connected component $H=\bid{C_n\circ K_p}$ with $n\geq 5$ and $p=\frac{1}{2}\bic(D)$, or every connected component $K$ of $D'$ admits an acyclic set $I_K$ (possibly empty) such that $\bic(K-I_K) \leq \bic(D') - 1$.
    
    In the former case, every vertex $x\in V(H)$ satisfies $d^+_H(x) = d^-_H(x) = \Delta$. By connectivity of $D$, some vertex $x\in V(H)$ is adjacent to $v$, implying that $\max(d^+_D(x),d^-_D(x)) > \Delta$, a contradiction to $\delmax(D) \leq \Delta$.
    In the latter case, taking the union of all such acyclic sets $I_K$ yields an acyclic set of $D'$ intersecting every maximum biclique of $D'$. Since $v$ does not belong to a maximum biclique of $D$, $I$ is also an acyclic set of $D$ such that $\bic(D-I) = \bic(D)-1$, a contradiction to $D$ being a counterexample. This shows that every vertex belongs to a maximum biclique.

    Let $\mathcal{H}$ be $\IG(\mathcal{C}(D))$, the intersection graph of the maximum bicliques of $D$. We denote by $\mathcal{H}_1,\dots,\mathcal{H}_r$ the connected components of $\mathcal{H}$.
    \begin{claim}
        \label{claim:empty_intersection}
        There exists a connected component $\mathcal{H}_i$ of $\mathcal{H}$ such that $\bigcap_{C\in V(\mathcal{H}_i)} C = \emptyset$.
    \end{claim}
    \begin{proofclaim}
        Assume for a contradiction that $X_i = \bigcap_{C\in V(\mathcal{H}_i)} C \neq \emptyset$ holds for every $i\in [r]$. Let $D'$ be the subdigraph of $D$ induced by $\bigcup_{i\in [r]} X_i$. 

        We claim that we can apply Theorem~\ref{thm:acyclic_hittingset} to $D'$, the vertices of $D'$ being partitioned into bicliques $X_1,\dots,X_r$, with $k = \left\lfloor \frac{1}{3}(\Delta+1)\right\rfloor$. 
        Let $i$ be any index of $[r]$ and $v_i$ be any vertex in $X_i$. Let us justify that the out-degree of $v_i$ outside $X_i$ in $D'$ is at most $k$, and that its in-degree outside $X_i$ is at most $|X_i| - k$.
        
        Observe that $v_i$ has at least $\bic(D)-1$ out-neighbours in $Y_i = \bigcup_{C\in V(\mathcal{H}_i)} C$, so the out-degree of $v_i$ outside $X_i$ in $D'$ is at most $\Delta - \bic(D)+1 \leq \frac{1}{3}(\Delta+1)$. Since it is an integer, the out-degree of $v_i$ outside $X_i$ in $D'$ is thus at most $k$. 
        On the other hand, $v_i$ has at most $\Delta - |Y_i|+1$ in-neighbours outside $X_i$ in $D'$. By Lemma~\ref{lemma:hajnal} applied to $S(D)$ we have $|Y_i| \geq 2\bic(D) - |X_i| \geq \frac{4}{3}(\Delta+1) - |X_i|$. Hence $v_i$ has at most $|X_i| - \frac{1}{3}(\Delta+1) \leq |X_i|-k$ in-neighbours outside $X_i$ in $D'$. 
        
        Therefore, we can apply Theorem~\ref{thm:acyclic_hittingset} to $D'$, and conclude that $D'$ admits an acyclic set $I$ of size $r$, which is thus a transversal of $(X_1,\dots,X_r)$. By definition of $D'$ and $(X_1,\dots,X_r)$, $I$ is thus an acyclic set of $D$ such that $\bic(D-I) = \bic(D)-1$, a contradiction.
    \end{proofclaim}

    Henceforth, by Claim~\ref{claim:empty_intersection}, we assume without loss of generality that $\bigcap_{C \in V(\mathcal{H}_r)}C = \emptyset$. By Lemma~\ref{lemma:christofides} applied to $S(D)$ with $\mathcal{C} = V(\mathcal{H}_r)$, we obtain $\bic(D) = \omega(S(D)) = \frac{2}{3}(\Delta+1)$ and $Y = \bigcup_{C\in \mathcal{C}} C$ admits a partition $Q_1,\dots,Q_n$ into $n\geq 4$ bicliques of $D$ each of cardinality $\frac{1}{2}\bic(D)$ such that $\mathcal{C} \setminus \{Q_1\cup Q_n\} = \{Q_i\cup Q_{i+1} : i\in [n-1] \}$. 
    
    If $Q_1\cup Q_n \in \mathcal{C}$, then $D\ind{Y}$ is $\Delta$-diregular, so necessarily $Y=V(D)$ as $D$ is connected, and $D$ is exactly $\bid{C_n\circ K_p}$. If $n$ is odd, this is a contradiction to $D$ being a counterexample. If $n$ is even, one can arbitrarily pick one representative of each $Q_i$ with $i\geq 1$ odd to obtain an acyclic set intersecting every maximum biclique of $D$, a contradiction.

    Hence $\mathcal{C} = \{Q_i\cup Q_{i+1} : i\in [n-1] \}$. Let $D'$ be the digraph obtained from $D$ by removing $\bigcup_{i=2}^{n-1} Q_i$ and adding every possible arc in both directions between $Q_1$ and $Q_n$. It is straightforward that $\delmax(D) \geq \delmax(D')$ and that $D'$ is smaller than $D$.

    Note that $\bic(D') \geq \bic(D)$ as $Q_1\cup Q_n$ is a biclique of $D'$ of size $\bic(D)$. We claim that equality holds, so assume for a contradiction that $D'$ contains a biclique $C'$ of size $\bic(D)+1$. 
    Let $q_1 = |C'\cap Q_1|$ and $q_n = |C' \cap Q_n|$. Observe that $q_1 > 0$ and $q_n >0$ as $C'$ is not a biclique of $D$.
    Let us fix $v\in C'\setminus (Q_1\cup Q_n)$, which exists since $|Q_1 \cup Q_n| < |C'|$. Recall that $v$ belongs to a maximum biclique $C$ of $D$, which is necessarily disjoint from $Q_1 \cup Q_n$ (otherwise $C$ would belong to $\mathcal{C}$).
    If $q_1+q_n > \frac{1}{3}(\Delta+1)$, then $d^+_D(v) \geq q_1+q_n+\bic(D)-1 > \Delta$, a contradiction. Hence $q_1+q_n \leq \frac{1}{3}(\Delta+1)$, which implies $|C'\setminus(Q_1\cup Q_n)| = \bic(D) +1 -q_1-q_n \geq \bic(D)+1 -\frac{1}{3}(\Delta+1)$. Let $v_1$ be any vertex in $C'\cap Q_1$, then 
    \[
    d^+_D(v_1) \geq |C'\setminus(Q_1\cup Q_n)| + |Q_1\cup Q_2|-1 \geq 2\bic(D) -\frac{1}{3}(\Delta+1) \geq \Delta+1,
    \]
    a contradiction.

    This shows that $\bic(D') = \bic(D)$. Note also that $D'$ is connected as $D$ is, and vertices in $\bigcup_{i=2}^{n-1}Q_i$ do not have neighbours outside $Y$. Since $D$ is a minimum counterexample, either $D' = \bid{C_{k}\circ K_p}$ or $D'$ has an acyclic set $I$ such that $\bic(D' - I) = \bic(D')-1$. 
    We distinguish the two possible cases.

    \medskip
    
    \noindent {\bf Case 1:} {\it The digraph $D'$ is isomorphic to $\bid{C_{k}\circ K_p}$.}

    \medskip

    In particular, observe that $D'$ is symmetric. Hence, by construction of $D'$, either $D$ is symmetric or it contains a simple arc $uv$ with $(u,v)\in Q_1\times Q_n$ or $(u,v)\in Q_n \times Q_1$. Indeed, the vertices in $Q_2\cup \ldots \cup Q_{n-1}$ are not incident to simple arcs, and any other simple arc of $D$ would be a simple arc of $D'$.
    
    If $D$ is symmetric, then the statement of the present theorem for $D$ corresponds exactly to the one of Theorem~\ref{thm:christofides} for $S(D)$. Hence, since $D$ is a counterexample, it cannot be symmetric, for otherwise $S(D)$ would be a counterexample to Theorem~\ref{thm:christofides}.
    
    Let thus $uv$ be a simple arc of $D$, and by symmetry assume that $u\in Q_1$ and $v\in Q_n$. 
    Since all arcs of $D'$ belong to $D$ except the ones between $Q_1$ and $Q_n$, we have
    \begin{align*}
        d^+_{D}(u) &= d^+_{D'}(u) - |Q_n\setminus N^+_D(u)| + |Q_2| &\\
        &\geq \Delta - (|Q_n|-1) +|Q_2| &\text{since $D' =\bid{C_{k}\circ K_p}$ and $v\in N^+_D(u) \cap Q_n$}\\
        &= \Delta+1,
    \end{align*}
    a contradiction to $\delmax(D)\leq \Delta$.

    \medskip
    
    \noindent {\bf Case 2:} \textit{The digraph $D'$ has an acyclic set $I$ such that $\bic(D' - I) = \bic(D')-1$.}

    \medskip

    Note that every maximum biclique of $D$ is either a maximum biclique of $D'$, or some set $Q_{i}\cup Q_{i+1}$ for $i\in [n-1]$.
    As $I$ intersects every maximum biclique of $D'$, in particular, it intersects $Q_1\cup Q_n$. We assume by symmetry that it intersects $Q_1$. For every $i\in [n]$, let $v_i$ be an arbitrary vertex of $Q_i$. If $2\leq i \leq n-1$, observe that $N^+(v_i) = N^-(v_i) = (Q_{i-1}\cup Q_i \cup Q_{i+1}) \setminus \{v_i\}$.

    If $n$ is even, $I \cup \{v_{2i+1} : i\in \left[ \frac{n-2}{2}\right]\}$ is an acyclic set intersecting every maximum biclique of $D$. 
    Note that it is indeed acyclic as $I\cap Q_n=\emptyset$, since $I$ intersects $Q_1$ and $Q_1\cup Q_n$ is a biclique of $D'$.
    
    If $n$ is odd, then $(I \setminus Q_1) \cup \{v_{2i} : i\in \left[\frac{n-1}{2}\right] \}$ is an acyclic set intersecting every maximum biclique of $D$. Note that this holds as vertices in $Q_1$ belong to exactly one maximum biclique in $D$, namely $Q_1\cup Q_2$. 
    Therefore, both cases yield a contradiction, and the result follows.
\end{proof}

\section{Handling vertices with a dense neighbourhood}
\label{sec:dense_neighbourhoods}

In the short proof of Theorem~\ref{thm:reed} due to King and Reed~\cite{kingJGT81}, Theorem~\ref{thm:sparse_graphs} handles the case of graphs in which every vertex has a sparse neighbourhood. The following result handles the remaining case, that is graphs admitting a vertex with a dense neighbourhood.

\begin{theorem}[{King and Reed~\cite[Thm.~6]{kingJGT81}}]
    \label{thm:dense_neighbourhoods}
    Let $a$ and $\epsilon$ be real numbers with $0<\epsilon < \frac{1}{6} - 2\sqrt{a}$. Let $G$ be a graph with maximum degree $\Delta $ satisfying $\omega(G) \leq \frac{2}{3}(\Delta+1)$ and let $v$ be such that $\edeg(v) > (1-a)\binom{\Delta}{2}$.
    Then 
    \[
    \chi(G) \leq \max\left( \chi(G-v), (1-\epsilon)(\Delta+1) \right).
    \]
\end{theorem}

The goal of this section is to prove Theorem~\ref{thm:handle_dense_vertex}, which is a directed analogue of Theorem~\ref{thm:dense_neighbourhoods} for digraphs. We first recall it here for convenience.

\thmdensevertex*

We use the following classical result due to Erd\H{o}s, Rubin, and Taylor~\cite{erdosCN26}.
\begin{lemma}[Erd\H{o}s {\it et al.}~\cite{erdosCN26}]
    \label{lemma:erdos_rubin_taylor}
    Let $G$ be a graph obtained from $K_n$ by removing a matching of size $p$. Then $G$ is $(n-p)$-choosable.
\end{lemma}

This directly gives the following. 

\begin{lemma}
    \label{lemma:generalisation_erdos_rubin_taylor}
    Let $D$ be a digraph obtained from $\bid{K_n}$ by removing a matching of size $p$. Then $D$ is $(n-p)$-dichoosable.
\end{lemma}
\begin{proof}
    Let $L$ be any $(n-p)$-list assignment of $D$. By Lemma~\ref{lemma:erdos_rubin_taylor}, there exists an $L$-colouring $\alpha$ of $S(D)$, because $S(D)$ is obtained from $K_n$ by removing a matching of size $p$. Since $D$ is obtained from $\bid{S(D)}$ by adding an acyclic set of arcs, $\alpha$ must be an $L$-dicolouring of $D$.
\end{proof}

\begin{proof}[Proof of Theorem~\ref{thm:handle_dense_vertex}]
    We take $\Delta(a) = \max\left(\frac{1-a}{\sqrt{a}-a},\frac{1-a}{a}\right)$.
    By directional duality we assume without loss of generality that $\edegp(v) > (1-a)\Delta(\Delta-1)$ and let $N=N^+(v) \cup \{v\}$ and $\bar{N} = V(D) \setminus N$. We assume that $|N| = \Delta+1$ for otherwise we add new vertices dominated by $v$.
    We partition $N$ into three parts as follows:
    \begin{align*}
        N_1 &= \left\{ u\in N : |N^+(u) \cap \bar{N}| \geq \frac{\Delta}{2} \right\},\\
        N_2 &= \bigg\{ u\in (N\setminus N_1) : |N^+(u) \cap (\bar{N}\cup N_1)| \geq 2\sqrt{a}\Delta \bigg\}, \text{and}\\
        N_3 &= N\setminus (N_1\cup N_2).
    \end{align*}
    Observe first that 
    \[\Delta^2\geq \sum_{u\in N\setminus\{v\}}d^+(u) \geq \edegp(v) + |A(N,\bar{N})|.\]
    Since $\edegp(v) > (1-a)\Delta(\Delta-1)$, we obtain 
    \[
        |A(N,\bar{N})| < \Delta^2 - (1-a)\Delta(\Delta-1) = a\Delta^2 + (1-a)\Delta.
    \]
    By definition of $N_1$, this implies
    \[
        \frac{\Delta}{2}|N_1| < a\Delta^2+(1-a)\Delta.
    \]
    Using the fact that $\Delta \geq \Delta(a) \geq \frac{1-a}{\sqrt{a}-a}$, we thus have $|N_1| < 2\sqrt{a}\Delta$.
    By definition, $v\notin N_1$ as $N^+(v) \subseteq N$. We also have $v\notin N_2$ as $N^+(v)\cap (\bar{N}\cup N_1) = N_1$, and $|N_1|< 2\sqrt{a}\Delta$.
    We now bound the size of $N_2$. Recall that every vertex $u \in N_1$ satisfies $2|N^+(u)\cap \bar{N}| \geq \Delta \geq d^-(u)$. We thus have the following.
    \begin{align*}
        |A(N_2\cup N_3, \bar{N}\cup N_1)| &= |A(N, \bar{N})| + |A(N_2\cup N_3, N_1)| - |A(N_1, \bar{N})| \\
        &\leq |A(N, \bar{N})| + \sum_{u\in N_1}d^-(u) - \sum_{u\in N_1}|N^+(u)\cap \bar{N}|\\
        &\leq |A(N, \bar{N})| + \sum_{u\in N_1}|N^+(u)\cap \bar{N}|\\
        &= |A(N, \bar{N})| + |A(N_1,\bar{N})|\\
        &\leq 2|A(N,\bar{N})| \\
        &< 2a\Delta^2 + 2(1-a)\Delta.
    \end{align*}
    By definition of $N_2$, this implies
    \[ 
    2\sqrt{a}\Delta|N_2| < 2a\Delta^2 + 2(1-a)\Delta,
    \]
    which together with $\Delta \geq \Delta(a) \geq \frac{1-a}{a}$ implies $|N_2|<2\sqrt{a}\Delta$.
    
    Let $k = \max(\dic(D-v), \left\lfloor (1-\epsilon)(\Delta+1)\right\rfloor)$. Since $v\in N_3$, $D-N_3$ is $k$-dicolourable, as it is a subdigraph of $D-v$. Let $\alpha$ be such a $k$-dicolouring and $L$ be the list assignment of $D\ind{N_3}$ defined as follows
    \[
    L(u) = [k] \setminus \alpha(N^+(u)\setminus N_3).
    \]    
    We now prove that $D\ind{N_3}$ is $L$-dicolourable, which implies that $D$ has dichromatic number at most $k$ (it is straightforward to check that any $L$-dicolouring of $D\ind{N_3}$ extends $\alpha$ into a dicolouring of $D$). First note that, by construction, every vertex $u\in N_3$ satisfies 
    \[
    |N^+(u) \setminus N_3| < 2\sqrt{a}\Delta +|N_2| < 4\sqrt{a}\Delta.
    \]
    Hence, every vertex $u\in N_3$ satisfies 
    \[
    |L(u)| > k -4\sqrt{a}(\Delta+1) \geq (1-\epsilon)(\Delta+1)-1-4\sqrt{a}(\Delta+1) \geq \frac{5}{6}(\Delta+1)-1.
    \]
    As $|L(u)|$ is an integer, we thus have $|L(u)| \geq \left\lfloor \frac{5}{6}(\Delta+1) \right\rfloor$.
    Let $M$ be a maximum matching in $\Bar{D}\ind{N_3}$ (the complement of $D\ind{N_3}$). Let $X$ be the vertices of $N_3$ that are not covered by $M$. By definition, $X$ is a biclique of $D$ so $|X|\leq \frac{2}{3}(\Delta+1)$.
    Note that $|N_3| = 2|M|+|X|$.  If $|M|\leq \frac{1}{6}(\Delta+1)$ then
    \[
    |M|+|X| \leq \frac{5}{6}(\Delta+1),
    \]
    and if $|M|\geq \frac{1}{6}(\Delta+1)$ then 
    \[
    |M|+|X| = |N_3|-|M| \leq \frac{5}{6}(\Delta+1).
    \]
    In both cases, we obtain $|M|+|X|\leq \frac{5}{6}(\Delta+1)$. Since $D\ind{N_3}$ is a subdigraph of a digraph obtained from $\bid{K_{|N_3|}}$ by removing a matching of size $|M|$, $D\ind{N_3}$ is $(|M|+|X|)$-dichoosable by Lemma~\ref{lemma:generalisation_erdos_rubin_taylor}. Therefore $D\ind{N_3}$ must be $\left\lfloor\frac{5}{6}(\Delta+1)\right\rfloor$-dichoosable, so in particular it is $L$-dicolourable, implying the result.
\end{proof}

\section{Proof of the main result}
\label{sec:proof_main}

This section is devoted to the proof of our main result, namely Theorem~\ref{thm:main}, which we first restate here for convenience.

\mainthm*

\begin{proof}
Let $\Delta_1$ be as in Corollary~\ref{cor:sparse_digraphs} and, for a real number $a >0$, let $\Delta(a)$ be as in Theorem~\ref{thm:handle_dense_vertex}. 
    We define $a$, $\epsilon_0$, $\Delta_2$, $\epsilon$, and $\gamma_{\epsilon}$ as follows:
    \begin{itemize}
        \item $a = \frac{1}{600}$,
        \item $\epsilon_0 = \frac{1}{6} - 4\sqrt{a} > 0$,
        \item $\Delta_2 = \min \left\{\Delta \geq 2 : \frac{\log^3 d}{d-1} <a \text{~for every integer $d\geq \Delta$}\right\}$,
        \item $\epsilon = \min\left( \frac{1}{\Delta_1},\frac{1}{\Delta_2}, \frac{1}{\Delta(a)},\frac{1}{3}\epsilon_0, \frac{a}{15e^7} \right)$, and
        \item $\gamma_\epsilon = \frac{\epsilon}{1-\epsilon} \leq 2\epsilon$.
    \end{itemize}
    We claim that the statement holds for this specific value of $\epsilon$. Assume that this is not the case, and let $D$ be a counterexample with minimum order.
    For the sake of conciseness let us denote $\deltil(D)$ by $\deltil$, $\delmax(D)$ by $\delmax$, $\bic(D)$ by $\bic$, and let $k=\lceil (1-\epsilon)(\deltil+1) + \epsilon \bic\rceil$. As $D$ is a minimum counterexample, it satisfies $\dic(D) > k$ and for every vertex $v\in V(D)$  we have $\dic(D-v) \leq k$. In particular, observe that $D$ is connected, for otherwise we can colour each connected component of $D$ independently with $k$ colours.
    
    This directly implies that every vertex $v$ of $D$ satisfies $\min(d^-(v),d^+(v)) \geq k$, for otherwise a $k$-dicolouring of $D$ can be obtained from any $k$-dicolouring of $D-v$ by choosing for $v$ a colour that is not appearing in its in-neighbourhood or in its out-neighbourhood.
    In particular, this implies that $\deltil$ must be close to $\delmax$, as we justify in the following claim.

    \begin{claim}
        \label{claim:deltil_close_delmax}
        $\delmax \leq (1+\gamma_\epsilon)\deltil$.
    \end{claim}
    \begin{proofclaim}
        Assume this is not the case, so there exists a vertex $v$ with $\max(d^+(v),d^-(v)) > (1+\gamma_\epsilon)\deltil$. By directional duality, we assume without loss of generality that $d^+(v) \geq d^-(v)$. Therefore, we have $\displaystyle \deltil^2 \geq d^+(v)\cdot d^-(v) > (1+\gamma_\epsilon)\deltil \cdot d^-(v)$, which implies 
        \[
        d^-(v) < \frac{1}{1+\gamma_\epsilon}\deltil = (1-\epsilon)\deltil \leq k,
        \]
        a contradiction.
    \end{proofclaim}
    
    With Theorem~\ref{thm:decrease_w} in hand, we can also justify that there exists a significant gap between $\bic$ and $\delmax$.

    \begin{claim}
        \label{claim:w_is_small}
        $\bic \leq \frac{2}{3}\left(\delmax +1\right)$.
    \end{claim}
    \begin{proofclaim}
        Assume for a contradiction that $\bic > \frac{2}{3}(\delmax+1)$. By Theorem~\ref{thm:decrease_w}, since $D$ is connected, there exists an acyclic set of vertices $I \subseteq V(D)$ which intersects every maximum biclique of $D$. Among all such sets $I$, we choose one with maximum size. In particular, this implies that every vertex in $V(D) \setminus I$ has at least one in-neighbour and one out-neighbour in $I$. Let $D'$ be $D-I$. As $I$ is acyclic and intersects every maximum biclique of $D$, we have $\bic(D') = \bic-1$. The following shows that we also have $\deltil(D') \leq \deltil -1$.
        \begin{align*}
            \deltil(D') &= \max \left\{ \sqrt{d^+_{D'}(v) \cdot d^-_{D'}(v)} : v\in V(D') \right\}\\
                    &\leq \max \left\{ \sqrt{(d^+_{D}(v)-1)(d^-_{D}(v)-1)} : v\in V(D') \right\}\\
                    &\leq \max \left\{ \sqrt{d^+_{D}(v)d^-_{D}(v)}-1 : v\in V(D') \right\}\\
                    &\leq \deltil-1,
        \end{align*}
        where in the first inequality we used the maximality of $I$ and in the second inequality we used that $\sqrt{(p-1)(q-1)} \leq \sqrt{pq}-1$ holds for every pair of real numbers $p,q\geq 1$. We briefly justify it here for completeness.
        
        \begin{equation*}
        \begin{alignedat}{2}
        & \quad & 0 &\leq (\sqrt{p} - \sqrt{q})^2\\
        \Ra && -p-q & \leq -2\sqrt{pq}\\
        \Ra && pq-p-q+1 & \leq pq-2\sqrt{pq}+1\\
        \Ra && (p-1)(q-1) &\leq (\sqrt{pq}-1)^2\\
        \Ra && \sqrt{(p-1)(q-1)} &\leq \sqrt{pq}-1 \\
        \end{alignedat}
        \end{equation*}
        Since $I$ is acyclic, we have $\dic(D) \leq \dic(D')+1$, as any dicolouring of $D'$ can be extended to $D$ by colouring $I$ with an additional colour. Together with the minimality of $D$ we obtain
        \begin{align*}
            \dic(D) &\leq \left\lceil(1-\epsilon)(\deltil(D')+1) + \epsilon\bic(D') \right\rceil + 1\\
            &\leq \left\lceil(1-\epsilon)\deltil + \epsilon(\bic-1) \right\rceil + 1\\
            &= \left\lceil(1-\epsilon)(\deltil+1) - (1-\epsilon) + \epsilon\bic -\epsilon + 1\right\rceil=k,
        \end{align*}
        a contradiction.
    \end{proofclaim}

    Our last claim guarantees that $\deltil$ is large enough.
    \begin{claim}
        \label{claim:large_delmax}
        $\delmax \geq \deltil \geq \max(\Delta_1,\Delta_2,\Delta(a))$.
    \end{claim}
    \begin{proofclaim}
        It is clear from the definitions that $\delmax \geq \deltil$.
        Assume for a contradiction that $\deltil < \max(\Delta_1,\Delta_2,\Delta(a))$.
        We assume that $\deltil >0$, the statement of the Theorem being trivial for $\deltil = 0$.
        Hence, since $\epsilon\leq \min\left(\frac{1}{\Delta_1}, \frac{1}{\Delta_2},\frac{1}{\Delta(a)}\right)$ we have $\epsilon < \frac{1}{\deltil}$, which implies:
            \[k = \left\lceil \deltil +1 + \epsilon(\bic - \deltil -1) \right\rceil \geq \left\lceil \deltil + \frac{1}{\deltil}(\bic-1)\right\rceil,\]  
        where in the inequality above we used that $\bic \leq \deltil + 1$. This holds as any vertex $v$ in a maximum biclique has in-degree and out-degree at least $\bic-1$, hence implying that 
        \[
        \deltil +1\geq \sqrt{(\bic -1)^2}+1 = \bic.
        \]
        Using the facts that $\dic(D) > k$ and that $\dic(D)\leq \left\lfloor\deltil + 1 \right\rfloor$, we deduce 
        \[\left\lfloor\deltil + 1 \right\rfloor > \left\lceil \deltil + \frac{1}{\deltil}(\bic-1)\right\rceil.\]
        This is a contradiction unless $\bic = 1$ and $\deltil$ is an integer. Henceforth we assume that $\bic=1$, \textit{i.e.}, $D$ is an oriented graph, and $\deltil \in \mathbb{N}$. 
        Therefore, since $\epsilon < \frac{1}{\deltil}$, we have
        \[
            k = \lceil(1-\epsilon)(\deltil +1) + \epsilon\rceil = \lceil \deltil + 1 - \epsilon \deltil \rceil = \deltil+1.
        \]
        But then $D$ satisfies
        \[
            \dic(D) \leq \deltil+1 = k,
        \]
        a contradiction. 
    \end{proofclaim}

    We now distinguish two cases, depending on the sparseness of $D$.

    \medskip
    
    \noindent{\bf Case 1:} \textit{There exists a vertex $v$ such that $\max(\edegm(v),\edegp(v)) > (1-a)\delmax(\delmax-1)$.}
    
    \medskip

    By Claim~\ref{claim:w_is_small} we have $\bic\leq \frac{2}{3}(\delmax+1)$. As $\epsilon_0 = \frac{1}{6}-4\sqrt{a}$, by Theorem~\ref{thm:handle_dense_vertex} we have 
    \[ \dic(D) \leq \max(\dic(D-v),(1-\epsilon_0)(\delmax+1) ). \]
    Recall that $\dic(D-v)\leq k$ as $D$ is a minimum counterexample. Since $\dic(D) > k$, we thus necessarily have 
    \[ \dic(D) \leq (1-\epsilon_0)(\delmax+1) \leq (1-\epsilon_0)(1+\gamma_\epsilon)(\deltil+1) \leq (1-\epsilon)(\deltil+1),\] 
    where in the second inequality we used Claim~\ref{claim:deltil_close_delmax}, and in the last inequality we used $\epsilon \leq \epsilon_0 - \gamma_\epsilon$, which holds because $\gamma_\epsilon \leq 2\epsilon$ and  $\epsilon \leq \frac{1}{3}\epsilon_0$. This yields the contradiction as $(1-\epsilon)(\deltil+1)\leq k$.

    \medskip
    
    \noindent{\bf Case 2:} \textit{For every vertex $v\in V(D)$ we have $\max(\edegm(v),\edegp(v)) \leq (1-a)\delmax(\delmax-1)$.}
    
    \medskip

    A fortiori every vertex satisfies $\min(\edegm(v),\edegp(v)) \leq (1-a)\delmax(\delmax-1)$.
    By Claim~\ref{claim:large_delmax} we know that $\delmax \geq \max(\Delta_1,\Delta_2)$. Hence, by definition of $\Delta_2$, we have $a>\frac{\log^3 \delmax}{\delmax-1}$. We can thus apply Corollary~\ref{cor:sparse_digraphs}, which together with Claim~\ref{claim:deltil_close_delmax} gives
    \[
    \dic(D) \leq \left(1-\frac{a}{5e^7}\right)(\delmax+1) \leq \left(1-\frac{a}{5e^7}\right)(1+\gamma_\epsilon)(\deltil+1) \leq \left(1-\epsilon\right)(\deltil+1),
    \]
    where in the last inequality we used $\epsilon \leq \frac{a}{5e^7}-\gamma_{\epsilon}$, which holds because $\gamma_\epsilon \leq 2\epsilon$ and $\epsilon \leq \frac{a}{15e^7}$. This contradicts $\dic(D)>k$ and concludes the proof. 
\end{proof}

\section{Further research: analogues for \texorpdfstring{$\delmin$}{Delta min}}
\label{sec:conclusion}

In this work we proved an analogue of Reed's result for digraphs by replacing the maximum degree of a graph with $\deltil$. A natural question then arises: can we strengthen this result by replacing $\deltil$ with $\delmin$? The answer to this question turns out to be negative. 
Moreover, for every fixed real numbers $\epsilon, \epsilon'$ with $0 < \epsilon' < 2\epsilon$, there exists a digraph $D$ with
\[
    \dic(D) > \lceil (1-\epsilon)(\delmin(D) + 1) + \epsilon'\bic(D) \rceil.
\]
This is because the {\sc $k$-Dicolourability} problem is NP-complete for every fixed $k\geq 2$, even when restricted to digraphs $D$ satisfying $\delmin(D) = k$ and $\bic(D) = \left\lceil \frac{k+1}{2} \right\rceil$, as shown in~\cite{picasarriJGT106}. In particular, the reduction shows that there exist such digraphs satisfying $\dic(D) = k+1$. Hence, if we choose $k$ to be at least $\frac{4-2\epsilon + 2\epsilon'}{2\epsilon-\epsilon'}$, we deduce the existence of digraphs $D$ with $\delmin(D) = k$ and $\dic(D) = k+1$, while $\left\lceil(1 - \varepsilon)(k + 1) + \varepsilon'\left\lceil\frac{k + 1}{2}\right\rceil\right\rceil \leqslant k$.

However, we believe that an analogue of our result holds with $\delmin$ instead of $\deltil$ if we also adapt the notion of clique number. Given a digraph $D$, a \textit{directed clique} of $D$ is a set of vertices $X$ that can be partitioned into $(X_1, X_2)$ such that both $D\ind{X_1}$ and $D\ind{X_2}$ are complete digraphs, and $D$ contains every possible arc $uv$ with $u\in X_1$ and $v\in X_2$. The \textit{directed clique number} of $D$, denoted by $\vec{\omega}(D)$, is the size of the largest directed clique of $D$.
We propose the following conjecture.

\begin{conjecture}
    \label{conj:delmin}
    There exists $\epsilon > 0$ such that every digraph $D$ satisfies 
    \[ \dic(D) \leq \lceil (1-\epsilon)(\delmin(D)+1) + \epsilon \vec{\omega}(D) \rceil. \]
\end{conjecture}

In particular, the conjecture above implies the following slightly weaker one, which, if true, shows that the condition $\epsilon' < 2\epsilon$ mentioned above is indeed necessary.

\begin{conjecture}
    There exists $\epsilon > 0$ such that every digraph $D$ satisfies \[ \dic(D) \leq \lceil (1-\epsilon)(\delmin(D)+1) + 2\epsilon \bic(D) \rceil. \]
\end{conjecture}

We are inclined to believe that these two conjectures are true even for $\epsilon = \frac{1}{2}$. One direction toward proving Conjecture~\ref{conj:delmin} would be to prove an analogue of Theorem~\ref{thm:main} using the maximum out-degree $\delplus$ instead of $\deltil$, that is to prove the existence of $\epsilon>0$ such that every digraph $D$ satisfies
\begin{equation}
    \label{eq:conj_delplus}
    \dic(D) \leq \lceil (1-\epsilon)(\delplus(D)+1) + \epsilon\bic(D) \rceil,
\end{equation}
where $\delplus(D) = \max_{v\in V(D)}(d^+(v))$. 
If such a result holds, then the following proposition, the proof of which uses a trick introduced in~\cite{picasarriJGT106}, implies Conjecture~\ref{conj:delmin}.

\begin{proposition}
    For every digraph $D$, there exists a digraph $H$ such that $\dic(D) \leq \dic(H)$, $\delplus(H)\leq \delmin(D)$ and $\bic(H) \leq \vec{\omega}(D)$.
\end{proposition}
\begin{proof}
    Let us fix a digraph $D$, and let $(X,Y)$ be the partition of $V(D)$ where $X=\{x\in V(D) : d^+_D(x) \leq \delmin(D)\}$ and $Y = V(D) \setminus X$. In particular, every vertex $y \in Y$ satisfies $d^-_D(y) \leq \delmin(D)$.

    Let $D'$ be the digraph obtained from $D$ by removing every arc from $Y$ to $X$ and replacing every arc from $X$ to $Y$ by a digon. Let $H$ be the digraph obtained from $D'$ by reversing every arc of $D'\ind{Y}$. It is straightforward to check that, by construction, $\delplus(H) \leq \delmin(D)$ and $\bic(H) \leq \vec{\omega}(D)$.

    We now justify that $\dic(H) \geq \dic(D)$. Assume that this is not the case, and let $\phi$ be a $(\dic(D) -1)$-dicolouring of $H$. By definition, $D$ coloured with $\phi$ must contain a monochromatic directed cycle $C = v_1,\dots,v_\ell,v_1$. Observe that $V(C) \nsubseteq X$, otherwise $C$ is a directed cycle of $H$. In a similar way, $V(C) \nsubseteq Y$, for otherwise $v_1,v_\ell,\dots,v_2,v_1$ is a directed cycle of $H$. Hence $C$ contains an arc $xy$ from $X$ to $Y$, but $\{xy,yx\}$ is a digon of $H$, a contradiction.
\end{proof}

We note that a result of Aharoni, Berger, and Kfir~\cite{aharoniJGT59} might be particularly useful to prove~\eqref{eq:conj_delplus}, as it provides a sufficient condition for a digraph with bounded maximum out-degree to admit an ASR (see~\cite[Corollary~II.13]{aharoniJGT59}). 
A natural direction would be to use this result to show that any minimum counterexample to~\eqref{eq:conj_delplus} must present a significant gap between $\delplus$ and $\bic$.
Nevertheless, the problem seems to be challenging already for oriented graphs. 

\begin{problem}
    Show the existence of $\epsilon>0$ such that every oriented graph $D$ with sufficiently large maximum out-degree $\delplus(D)$ satisfies
    $\dic(D) \leq (1-\epsilon)\delplus(D)$.
\end{problem}

To prove the analogue of the problem above with $\deltil$ instead of $\delplus$, two distinct strategies are mainly known. As we explain below, neither of them works when only the out-degrees are bounded; the problem above thus requires new methods.

The first one, which we used in Section~\ref{sec:oriented_graphs}, consists of decomposing the given digraph into digraphs with small maximum degree (thanks to Lov\'asz' result, \Cref{thm:lovasz}) and then applying an analogue of Brooks' theorem on each smaller digraph. Unfortunately, the only directed analogue of \Cref{thm:lovasz} involving $\delplus$, due to Alon~\cite{alonCPC15} (see also~\cite{alonJCTB68}), requires an extra factor $2$, which prevents this strategy from working. 

The second one, which we used in Section~\ref{sec:sparse_digraphs}, is the probabilistic method, involving the Lov\'asz Local Lemma.
This is doomed to fail, as if the in-degree can be arbitrarily larger than the out-degree, one cannot expect to bound the degree of ``bad events'' in order to apply Lov\'asz    Local Lemma.
This drawback lies at the heart of some long-standing conjectures in digraph theory; the most famous one is probably the conjecture of Stiebitz~\cite{stiebitzKAM95} (popularized by Alon~\cite{alonJCTB68}) that if a digraph $D$ has sufficiently large minimum out-degree, then it can be vertex-partitioned into two digraphs with minimum out-degree at least $2$.

\section*{Acknowledgments}
Research supported by JSPS Kakenhi 26K21777 and JP25K24465 and by JST ASPIRE JPMJAP2302.

\bibliographystyle{abbrv}
\bibliography{refs}

\end{document}